\documentclass[leqno,11pt]{amsart}
\usepackage{latexsym,amsfonts,amsmath}
\usepackage{mathrsfs}
\usepackage{graphicx,graphics,epsf,psfrag,eepic,epic}
\usepackage{xypic}
\xyoption{all}

\input epsf.tex

\theoremstyle{plain}
\newtheorem{theorem}{Theorem}[section]
\newtheorem{corollary}[theorem]{Corollary}
\newtheorem{proposition}[theorem]{Proposition}
\newtheorem{lemma}[theorem]{Lemma}

\theoremstyle{definition}
\newtheorem{definition}[theorem]{Definition}
\newtheorem{remark}[theorem]{Remark}
\newtheorem{example}[theorem]{Example}

\newcommand{\C}{{\mathbb C}}
\newcommand{\R}{{\mathbb R}}
\newcommand{\Z}{{\mathbb Z}}
\newcommand{\CP}{{\mathcal{P}}}
\newcommand{\CR}{{\mathcal{R}}}
\newcommand{\CA}{{\mathcal{A}}}
\newcommand{\CF}{{\mathcal{F}}}
\newcommand{\CK}{{\mathcal{K}}}
\newcommand{\CH}{{\mathcal{H}}}
\newcommand{\CI}{{\mathcal{I}}}
\newcommand{\CB}{{\mathcal{B}}}
\newcommand{\CC}{{\mathcal{C}}}
\newcommand{\CL}{{\mathcal{L}}}

\newcommand{\CG}{{\mathcal{G}}}
\newcommand{\CD}{{\mathcal{D}}}
\newcommand{\CV}{{\mathcal{V}}}

\newcommand{\gc}{{\mathfrak{c}}}
\newcommand{\vol}{{\mathrm{vol}}}
\newcommand{\res}{{\mathrm{res}}}
\newcommand{\Ires}{{\mathrm{Ires}}}
\renewcommand{\ll}{{\langle}}
\newcommand{\rr}{{\rangle}}

\renewcommand{\c}{\mathfrak{c}}
\def\spa{\operatorname{Sp}(a)}
\def\h{\operatorname{ht}}
\def\etc{\emph{etc.} }
\def\lra{\longrightarrow}
\def\jac{\operatorname{Jac}}
\def\Hyp{\operatorname{Hyp}}
\def\JK{\operatorname{JK}}
\def\maxi{\mathit{maxi}}
\def\la{\lambda}
\def\vector#1#2{\left(\begin{array}{cc}#1\\#2\\ \end{array}\right)}
\def\rt{z_{i_0}-z_{j_0}}
\def\binome#1#2{{\genfrac{(}{)}{0pt}{0}{#1}{#2}}}
\def\bases{\operatorname{Bases}}
\def\sec{\operatorname{s}}
\def\mo{\operatorname{Mo}}
\def\proj{\operatorname{proj}}

\title[Partition functions for classical root systems]{Volume
  computation for polytopes and partition functions for
  classical root systems}
\author[M.W. Baldoni, M. Beck, C. Cochet, M. Vergne]{M.\ Welleda Baldoni, Matthias Beck, Charles Cochet, Mich\`ele Vergne}

\date{June 15, 2005}


\thanks{During the course of this work, we have benefitted from
  discussions with Jes\'us de Loera, Andr\'as Szenes, Corrado de
  Concini and Claudio Procesi. We would like to thank them for sharing
  their mathematical expertise with us. 
We also thank the  various institutions that helped us to collaborate
  on this work:\ the Research-in-pairs program at the
  Forschungsinstitut Oberwolfach, a LIEGRITS grant and the University
  Tor Vergata Roma, the University Denis Diderot in Paris, the Centre
  Laurent Schwartz at Ecole Polytechnique, and San Francisco State
  University.
}

\begin{document}

\maketitle

\begin{abstract}
This paper presents an algorithm to compute the value of the
inverse Laplace transforms of rational functions with poles on
arrangements of hyperplanes. As an application, we present an
efficient  computation of the partition function for classical
root systems.
\end{abstract}

\section{Introduction} \label{sect.intr}

The ultimate goal of this work is to present an algorithm for a fast
computation of the partition function of classical root systems.
We achieve this goal in somewhat more general terms, namely we
develop algorithms to compute the volume of a polytope and its
discrete analog, the number of integer points in the
polytope. These formulas, in turn, are inverse Laplace
transforms of certain rational functions, and our work can be viewed in
these general terms.

Let $U$ be a finite-dimensional real vector space of dimension $r$.
Denote its dual vector space $U^*$ by $V$.
Consider a set of elements
$$
\CA=\{\alpha_1,\alpha_2,\ldots,\alpha_N\}
$$
of non-zero vectors of $V$.
We assume that the convex cone $\CC(\CA)$ generated by non-negative
linear combinations of the elements $\alpha_i$ is an acute convex
cone in $V$ with non-empty interior.

The elements $\ell$ in $V$ produce linear functions
$u\mapsto\ell(u)$ on the complexified vector space $U_\C$.
In particular, to the set $\CA$ we associate the arrangement of
hyperplanes
$$
\CH_{\C}(\CA):=\bigcup_{i=1}^N\left\{u\in U_\C\,|\,\,\alpha_i(u)=0\right\}
$$
in $U_\C$ and its complement
$$
U_{\C}(\CA):=\left\{u\in U_\C\,\Big|\,\,
    \prod_{i=1}^N\alpha_i(u)\neq 0\right\}.
$$
We denote by $\CR_{\CA}$ the ring of rational functions on $U_\C(\CA)$
with poles along $\CH_{\C}(\CA)$.
Then each element $\phi\in\CR_{\CA}$ can be written as $P/Q$
where $P$ is a polynomial function on $r$ complex variables and
$Q$ is a product of elements, not necessarily distinct, of $\CA$.

Our first aim is to present an algorithm to compute the value
of the inverse Laplace transform of functions in $\CR_{\CA}$ at a
point $h\in V$. In other words, we study the value at a point $h\in
V$ of convolutions of a number of Heaviside distributions
$\phi\mapsto\int_0^{\infty}\phi(t\alpha_i) dt$. The first
theoretical ingredient is the notion of Jeffrey-Kirwan
residues~\cite{jeffreykirwan}. Going a step further,
DeConcini-Procesi~\cite{deconciniprocesi} proved that one can
compute Jeffrey-Kirwan residues using maximal nested sets (in short
MNS), a combinatorial tool related to no-broken-circuit bases of
the set of vectors $\CA$.

The applications in view are volume computation for polytopes,
enumeration of integral points in polytopes and, more generally,
discrete or continuous integration of polynomial functions over
polytopes. Indeed, Szenes-Vergne~\cite{SzeVer}, refining a formula
of Brion-Vergne~\cite{BriVer97}, stated formulae for the volume
and number of integral points in polytopes involving
Jeffrey-Kirwan residues.

 Consider the polytope
$$\Pi_{\CA}(h)
  :=\left\{x\in\R^N
          \,\Big|\,
          \sum_{i=1}^Nx_i\,\alpha_i=h,\,x_i\geq 0\right\}.
$$
 As a function of $h$, the volume of $\Pi_{\CA}(h)$ is a
piecewise-defined polynomial. The chambers of polynomiality in the
parameter space $V$ are polyhedral cones.

Our programs are extremely efficient for computing the volume
 of the polytope $\Pi_{\CA}(h)$ when $\CA$ is a classical root
 system.
An important fact is that our algorithm can work with formal
parameters, thus giving the polynomial volume formula for
$\Pi_{\CA}(h)$ when $h$ runs over a particular chamber.

For an analogous theory for integral-point enumeration, we have to
assume that the $\alpha_i$ are  vectors in a lattice $V_\Z$. For
$h\in V_\Z$, the function $N_{\CA}(h)$ which associates to the
vector $h$ the number of integral points in $\Pi_{\CA}(h)$, that
is the number of ways to represents the vector $h$ as a sum of a
certain number of vectors $\alpha_i$, is called the
(vector)-partition function of $\CA$. For example for $B_2$, given
a vector $(h_1,h_2)$ with integral coordinates we
would like to compute the number $N_{B_2}(h)$ of vectors
$(x_i)\in\Z_+^4$ such that
$$x_1\vector{1}{0}+x_2\vector{0}{1}
 +x_3\vector{1}{-1}+x_4\vector{1}{1}=\vector{h_1}{h_2}.$$

As a function of $h$, the number $N_{\CA}(h)$ of integral points
in $\Pi_{\CA}(h)$ is a piecewise-defined quasipolynomial, and
again the chambers of quasipolynomiality are polyhedra in
$V$~\cite{sturmfels,Sze}.

In this paper, we describe an efficient algorithm for MNS
computation for classical root systems. This algorithm for MNS
gives rise to programs for Kostant partition function for the
classical root systems $A_n$, $B_n$, $C_n$, and $D_n$. Again, our
algorithm works with a formal parameter $h$ that is assumed to be
confined to a particular chamber.

These calculations are valuable because partition functions play a
fundamental role also in representation theory of semisimple Lie
algebras $\frak g$. Indeed, partition functions arise naturally
when we want to compute the multiplicity of a weight in a
finite-dimensional representation or the tensor-product
decomposition of two representations, both being basic problems to
understand characters of representations. Cochet~\cite{cochet} has
obtained very efficient algorithms for both these problems in the
case of $A_n$, implementing results of~\cite{BalDeLoeVer03}. See
also a forthcoming paper \cite{cochetABCD} for
multiplicities computation in all the classical Lie algebras using
the results obtained in this paper. There is also a class of
infinite-dimensional representations, the discrete-series
representations, whose understanding  is central for the general
theory of admissible irreducible representations. The
decomposition of such representations to a maximal compact
subgroup of $\frak g$ is predicted by Blattner's formula, which is
a partition function in which the roots involved are the so-called
noncompact roots.

We conclude by describing the   way  the paper is organized.
 Section~\ref{sect.lapl}
introduces Laplace transforms and polytopes. In
Section~\ref{sect.def.JK}, we recall Jeffrey-Kirwan residues and
its link with counting formulae.
DeConcini-Procesi's maximal nested sets are described
in Section~\ref{sect.JK}, as well as how they are related to
Jeffrey-Kirwan residues. Section~\ref{sect.gene} describes our
general algorithm for MNS computations. Details of particular
cases of the algorithm for the root systems  $A_n$, $B_n$, $C_n$
and $D_n$ are examined in Sections~\ref{sect.A}--\ref{sect.D}.
Finally comparative tests of our programs with existing softwares
are performed in Section~\ref{sect.perf}.

A number of  theoretical results  on the function $N_{\CA}(h)$
when $\CA$ is a subset of the system $A_n$  can be found  in
Baldoni-Vergne \cite {baldonivergne} (as, for example, the
computation of the volume of the Chan-Robbins polytope).

Computer programs for volume computation/integral-point enumeration in
polytopes have only been implemented in the very recent past, most 
notably {\tt LattE}~\cite{lattemanual,latte} and {\tt
barvinok}~\cite{BBLSV04}, both of which are implementations of
Barvinok's algorithm~\cite{barvinok}. To the best of our
knowledge, these two are the only general programs for volume
computation/integral-point enumeration in polytopes. More
specialized programs include algorithms of Baldoni-DeLoera-Vergne
for flow polytopes~\cite{BalDeLoeVer03} and Beck-Pixton for the
Birkhoff polytope~\cite{beckpixton}.

Our programs have been especially designed for classical root
systems, are faster than all actual existing softwares and can compute
new examples that were not reachable by previous algorithms.
Note in particular that our programs can perform computations for
$N_{\CA}(h)$ for $\CA_n$ at least up to $n=10$ (11 coordinates
vector).
For $\CB_n$, $\CC_n$, $\CD_n$ the algorithms are efficient at least up
to $n=6$.
For our methods (as well as for {\tt LattE}), the size of the vector
$h$ affects only little on the computation time.
Recall that our methods can also calculate the multivariate
quasi-polynomials $h\mapsto N_{\CA}(h)$ when $h$ varies on a chamber,
and as a particular case for a  fixed $h$ the function 
$k\mapsto N_{\CA}(k h)$ which is the Ehrhart quasipolynomial in $k$.

\section{Laplace transform and polytopes} \label{sect.lapl}

We start by briefly recalling the notations of the introduction,
aiming to relate the  inverse of the Laplace transform with
various \emph{counting formulae} for a polytope. A good
introduction on this theme is the survey article~\cite{ver}.

\subsection{Laplace transform} \label{sect.defi}

Let $U$ be a finite-dimensional real vector space of dimension $r$
with dual space $V$. We fix the choice of a Lebesgue measure $dh$
on $V$.
Consider a set
$$\CA=\{\alpha_1,\alpha_2,\ldots,\alpha_n\}$$
of non-zero vectors of $V$. We assume that the set of vectors
$\alpha_i$ spans $V$. For any subset $S$ of $V$, we denote by
$\CC(S)$ the convex cone generated by non-negative linear
combinations of elements of $S$. We assume that the convex cone
$\CC(\CA)$ is acute in $V$ with non-empty interior.

\begin{figure}[ht]
\begin{center}
\includegraphics[height=4cm,width=6cm]{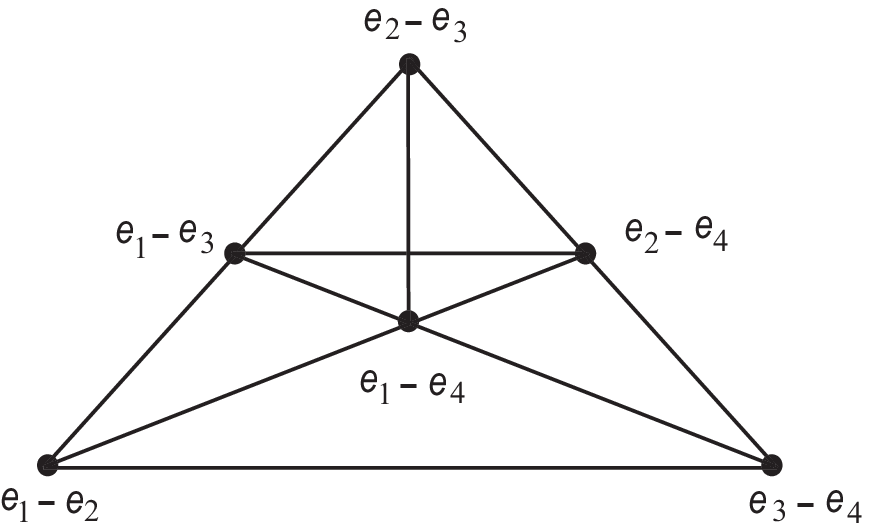}
\caption{The 7 chambers for $A_3$}
\label{figu.A3.cham}
\end{center}
\end{figure}

\begin{figure}[ht]
\begin{center}
\includegraphics[height=5cm,width=7.5cm]{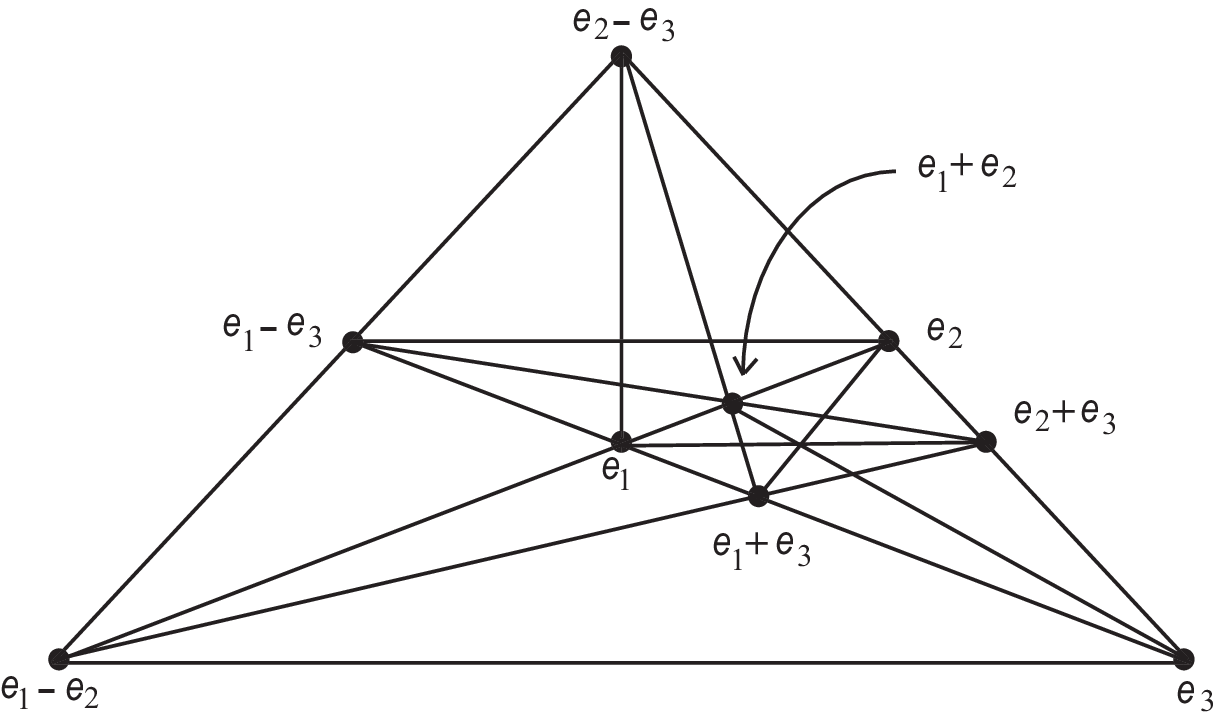}
\caption{The 23 chambers for $B_3$} \label{figu.B3.cham}
\end{center}
\end{figure}

Let $\CV_{sing}(\CA)$ be the union of the boundaries of the cones
$\CC(S)$, where $S$ ranges over all the subsets of $\CA$. The
complement of $\CV_{sing}(\CA)$ in $\CC(\CA)$ is by definition the
open set $\CC_{reg}(\CA)$ of \emph{regular} elements. A connected
component $\gc$ of $\CC_{reg}(\CA)$ is called a \emph{chamber} of
$\CC(\CA)$. Figures~\ref{figu.A3.cham} and~\ref{figu.B3.cham}
represent slices of the cones $\CC(A_3)$ and $\CC(B_3)$, where the
dots represent the intersection of a slice with a ray $\R_{ \ge 0
}$ hence showing the chambers. Note that the chambers for $B_r$
and $C_r$ are the same (as roots in  $B_r$ and $C_r$ are
proportional). In dimension $3$, the root system  $A_3$ is
isomorphic to $D_3$. See~\cite{BalDeLoeVer03} for the computation
of chambers. Very little is known about the total number of
chambers. On the other hand, given a vector $h$, it is easy to
compute the equations of the chamber containing $h$. This was done
in~\cite{BalDeLoeVer03,cochethome}. We have incorporated this small part of
the corresponding program in our programs for classical root
systems.

Table~\ref{figu.chambers} represents the only numbers of chambers
that have been computed (and the computation time).

\begin{figure}[ht]
\begin{tabular}{|c|c|c|c|c|c|c|} \hline
   &      A   &      B   &    C     &    D
&    F    &    G    \\ \hline
 1 &         1&         1&       1  &
&         &         \\
   &      (0s)&      (0s)&    (0s)  &
&         &         \\ \hline
 2 &         2&         3&         3&       1
&         &       5 \\
   &      (0s)&      (0s)&      (0s)&    (0s)
&         &    (0s) \\ \hline
 3 &         7&        23&        23&       7
&         &         \\
   &      (1s)&      (8s)&      (8s)&    (1s)
&         &         \\ \hline
 4 &        48&       695&       695&     133
&  12946  &         \\
   &     (23s)&     (11m)&     (11m)&   (90s)
& (3d16h) &         \\ \hline
 5 &       820&  $>$26905&  $>$26905&   12926
&         &         \\
   &     (19m)&  ?&  ?&  (1d5h)
&         &         \\ \hline
 6 &     44288&         ?&         ?&       ?
&         &         \\
   &  (24d18h)&          &          &
&         &         \\ \hline
\end{tabular}
\caption{Number of chambers and computation time}
\label{figu.chambers}
\end{figure}

Consider now a cone $\CC(S)$ spanned by a subset $S$ of $\CA$ and
let $p$ be a function on $\CC(S)$. We assume that $p$ is the
restriction to $\CC(S)$ of a polynomial function on $V$. By
superposing such functions $p$, we
obtain a space $\CL\CP(V,\CA)$ of locally polynomial functions on
$\CC(\CA)$. For $f\in\CL\CP(V,\CA)$, the restriction of $f$ to any
chamber $\gc$ of $\CC(\CA)$ is given by a polynomial function.

The Laplace transform $L(f)$ of such a function $f$ is defined as
follows.
Consider the dual cone $\CC(\CA)^*\subset U$ of $\CC(\CA)$ defined
by:
$$\CC(\CA)^*=\left\{u\in U\,|\,\,\ll h,u\rr\geq 0
  \mbox{ for all }h\in\CC(\CA)\right\}.$$
Then for $u$ in the interior of the cone $\CC(\CA)^*$, the
integral
$$L(f)(u)=\int_{\CC(\CA)}e^{-\ll h,u\rr}f(h)dh$$
is convergent.
It is easy to see that the function $L(f)$ is the restriction
to $\CC(\CA)^*$ of a function in $\CR_{\CA}$.
(Recall that $\CR_{\CA}$ is the ring of rational
functions $P/Q$ on $U$ where $P$ is a polynomial function on $U$
and  $Q$ is a product of elements of $\CA$.)
It is easy~\cite{BriVer97} to characterize the functions
$L(f)$ on $U$ arising this way.

Let $\nu$ be a subset of $\{1,2,\ldots,n\}$.
We will say that $\nu$ is \emph{generating} (respectively
\emph{basic}) if the set $\{\alpha_i\,|\,i\in\nu\}$ generates
(respectively is a basis of) the vector space $V$.

Every basic subset is of cardinality $r$ and we write
$\bases(\CA)$ for the set of basic subsets.
Given $\sigma\in\bases(\CA)$, the associated \emph{basic
fraction} is
\begin{equation} \label{equa.basi}
f_\sigma=\frac{1}{\prod_{i\in\sigma}\alpha_i}.
\end{equation}
In a system of coordinates (depending on $\sigma$) on $U$ where
$\alpha_i(u)=u_i$ (for $i\in\sigma$), such a basic fraction is
simply of the form
$$\frac{1}{u_1u_2\cdots u_r}.$$

Define $\CG(U,\CA)\subset\CR_{\CA}$ as the linear span of
functions $\frac{1}{\prod_{i\in\nu}\alpha_i^{n_i}}$, where $\nu$
is generating and $n_i$ are positive integers.
The following proposition gives the characterization we were
speaking of and is easy to prove:

\begin{proposition}\cite{BriVer97} \label{prop.bv}
If $f$ is a locally polynomial function on $\CC(\CA)$, the Laplace
transform $L(f)$ of $f$ is the restriction to $\CC(\CA)^*$ of a
function in $\CG(U,\CA)$.
Reciprocally, for any generating set $\nu$ and every set of
positive integers $n_i>0$, there exists a locally polynomial
function $f$ on $V$ such that
$$\frac{1}{\prod_{i\in\nu}\alpha_i(u)^{n_i}}
  =\int_{\CC(\CA)}e^{-\ll h,u\rr}f(h)dh$$
for any $u$ in the interior of $\CC(\CA)^*$.
\end{proposition}

We define the inverse Laplace transform
$L^{-1}:\CG(U,\CA)\to\CL\CP(V,\CA)$ as follows.
For $\phi\in\CG(U,\CA)$, the function $L^{-1}\phi$ is the
unique locally polynomial function that satisfies
$$\phi(u)=
  \int_{\CC(\CA)}e^{-\ll h,u\rr}(L^{-1}\phi)(h)dh$$
for any $u\in\CC(\CA)^*$.

In the next sections, we will explain the relation between Laplace
transforms and the enumeration of integral points of families of
polytopes. We will see in Section~\ref{sect.JK} that one can write
efficient formulae for the inversion of Laplace transforms in
terms of residues, whose algorithmic implementation is working in
a quite impressive way, at least for low dimension.

\subsection{Volume and number of integral points of a polytope}
\label{sect.volu}

In this subsection we consider a sequence
$$\CA^+=[\alpha_1,\alpha_2,\ldots,\alpha_N]$$
of non-zero elements of $\CA$. We assume that each element
$\alpha\in\CA$ occurs in the sequence; in particular $N\geq n$ and
the set $\CA^+$ spans $V$.

\begin{remark} \label{rema.mult} In all our examples, the sequence
$\CA^+$ will not have multiplicities, so that we will freely
identify $\CA^+$ and $\CA$.
\end{remark}

 We introduce now the notion of a
partition polytope.

We consider the space $\R^N$ with its standard basis $\omega_i$
and Lebesgue measure $dx$.

If $x=\sum_{i=1}^N x_i\omega_i\in\R^N$ with $x_i\geq 0$
($1\leq i\leq N$) then we will simply write $x\geq 0$.

Consider the surjective map $A:\R^N\to V$ defined by
$A(\omega_i)=\alpha_i$ and denote by $K$ its kernel.
Then $K$ is a vector space of dimension $d=N-r$ equipped
with the quotient Lebesgue measure $dx/dh$.

If $h\in V$, we define
$$\Pi_{\CA^+}(h)=\left\{x\in\R^N\,|\,Ax=h;x\geq 0\right\}.$$
The set $\Pi_{\CA^+}(h)$ is a convex polytope.
It is the intersection of the non-negative quadrant in $\R^N$ with
an affine translate of the vector space $K$.
This polytope consists of all non-negative solutions of the system
of $r$ linear equations
$$\sum_{i=1}^N x_i\alpha_i=h.$$

\begin{remark} \label{rema.poly}
It might be appropriate to recall that any full dimensional convex
polytope $P$ in a vector space $E$ of dimension $d$, defined
by a system of $N$ linear inequations
$$P=\{y\in E\,|\,\ll u_i,y\rr+\lambda_i\geq 0\}$$
(where $u_i\in E^*$ and $\lambda_i$ are real numbers), can be
canonically realized as a partition polytope $\Pi_{\CA^+}(h)$.
Here $\CA^+$ is a sequence of $N$ elements in a vector space of
dimension $r=N-d$. Indeed, consider the diagram
$$E\stackrel{i}{\lra}\R^N\stackrel{A}{\lra}V=\R^N/i(E)$$
where $i:y\mapsto\sum_{i=1}^N\ll u_i,y\rr\omega_i$ and
$A$ is the projection map $\R^N\lra V$. Let $\alpha_i$ be the
images of the canonical basis $\omega_i$ of $\R^N$. Define
$\CA^+=[\alpha_1,\ldots,\alpha_N]$ and consider the point
$h:=A(\sum_{i=1}^N\lambda_i\omega_i)$. Then the polytope
$\Pi_{\CA^+}(h)$ is isomorphic to $P$. Indeed, the points in
$\Pi_{\CA^+}(h)$ are exactly the points $x_i$ such that
$\sum_{i=1}^N (x_i-\lambda_i) A(\omega_i)=0$ with $x_i\geq 0$. By
definition of the space $V=R^N/i(E)$, there exists $y\in E$ such
that $x_i-\la_i=\ll u_i,y\rr$. As $x_i\geq 0$, this means exactly
that $\ll u_i,y\rr+\la_i\geq 0$, so that the point $y$ is in $P$.

More concretely, to determine the partition polytope $Ax=b$
starting from a polytope $P$ given by $Qy^T\geq\la$
(where $Q$ is a $N\times d$ matrix whose $i^{\text{th}}$ row is
given by a vector $u_i\in E^*$ and $\la\in\R^N$) we choose
among the elements $u_i$ a basis of $E^*$.
Thus after relabeling the indices and doing an appropriate
translation, we may assume the inequations of the polytope $P$
are given in the form
$$\begin{array}{lll}
\left\{\begin{array}{l}
  y_1\geq 0\\
  y_2\geq 0\\
  \ldots\\
  y_d\geq 0\\
  C y^T+\la^T\geq 0
\end{array}\right.
&\text{where}&
  \begin{array}{ll}
    \mbox{$C$ is a $r\times d$ matrix,}\\
    \la\in\R^r,\\
    \mbox{and $y=(y_1,\ldots,y_d)\in\R^d$.}\\
  \end{array}
\end{array}$$
Then the polytope $P$ is isomorphic to the polytope defined by
$$\left\{x\geq 0\,|\,\,Ax^{\mbox{T}}=\la^{\mbox{T}}\right\}$$
where  $A$ is the $r\times N$ matrix given by
$$A=\begin{array}{cr}
\left(
\begin{array}{cc}{}
\underbrace{-C}_{r\times d}&\underbrace{I_r}_{ r\times r}\
\end{array}
\right ),
& I_r\ \ \text{being the identity matrix.}
\end{array}$$
\end{remark}

\begin{example} \label{exam.poly}
Let $P\subset\R^2$ be the polytope defined by the system
of inequalities:
$$\left\{\begin{array}{rcl}
-x_1+1&\geq&0,\\
-x_2+2&\geq&0,\\
-x_1-x_2+2&\geq&0,\\
2x_1+x_2-1&\geq&0.\\
\end{array}\right.$$
Choosing the basis $u_1=(1,0)$, $u_2=(0,1)$ and using the
translation $y_1=-x_1+1$ and $y_2=-x_2+2$ we can rewrite the
system as:
$$\begin{array}{lr}
  \left\{\begin{array}{l}
    y_1\geq 0\\
    y_2\geq 0\\
    C\left(\begin{array}{c}y_1\\ y_2\\ \end{array}\right )
    +\left(\begin{array}{c}-1\\ 3\end{array}\right)\geq 0
  \end{array}\right.
&
\mbox{where}\ \
C=\left(\begin{array}{cc}1&1\\ -2&-1\end{array}\right).
\end{array}$$
Therefore $P$ is isomorphic to
$$\Pi_{\CA^+}(h)=\Big\{y=(y_1,\ldots,y_4)\in\R^4,y\geq 0
\ \Big | \
\begin{array}{rcl}-y_1-y_2+y_3&=&-1\\ 2y_1+y_2+y_4&=&3\\ \end{array}
\Big\}
$$
with $h=\vector{-1}{3}$.
\end{example}

We continue with our review.
If $h$ is in the interior of the cone $\CC(\CA)$, then the
polytope $\Pi_{\CA^+}(h)$ is of dimension $d$.
It lies in a translate of the vector space $K$, and this
translated space is provided with the quotient measure $dx/dh$.

\begin{definition} \label{defi.volu}
We write $\vol_{\CA^+}(h)$ for the volume of $\Pi_{\CA^+}(h)$
computed with respect to this measure.
\end{definition}

Suppose further that $V$ is provided with a lattice
$V_\Z$ and that
$$\CA^+:=[\alpha_1,\alpha_2,\ldots,\alpha_N]$$
is a sequence of non-zero elements of $V_\Z$ spanning
$V_\Z$, that is, $V_\Z = \sum_{ i=1 }^N \Z \alpha_i$.

In this case, the lattice $V_\Z$ determines a measure $d_{\Z}h$ on
$V$ so that the fundamental domain of the lattice $V_\Z$ is of
measure $1$ for $d_{\Z}h$. However, for reasons which will be
clear later on, we keep our initial  measure $dh$. We introduce
the normalized  volume.

\begin{definition} \label{defi.voluZ}
The normalized volume  $\vol_{\Z,\CA^+}(h)$ is the volume of
$\Pi_{\CA^+}(h)$ computed with respect to the measure
$dx/d_{\Z}h$.
\end{definition}

\begin{remark}
The reason for keeping our initial $dh$ is that the root systems
$B_r, C_r,D_r$ live on the same standard vector space $V=\R^r$,
where the most natural measure is the standard one. This measure
is twice the measure given by the root lattice in the case of
$C_r$ and $D_r$.
\end{remark}

If $\vol(V/V_\Z,dh)$ is the volume of a fundamental domain of
$V_\Z$ for $dh$, clearly $
\vol_{\Z,\CA^+}(h)=\vol(V/V_\Z,dh)
\vol_{\CA^+}(h).$

Let now $h\in V_\Z$. A discrete analogue of the normalized volume
of $\Pi_{\CA^+}(h)$  is the number of integral points inside this
polytope.

\begin{definition} \label{defi.numb}
Let $N_{\CA^+}(h)$ be the number
of integral points in $\Pi_{\CA^+}(h)$, that is the number
of solutions $x=(x_1,\ldots,x_N)$ of the equation
$\sum_{i=1}^N x_i\alpha_i=h$ where $x_i$ are non-negative
integers.
The function $h\mapsto N_{\CA^+}(h)$ is called the partition function
of $\CA^+$.
\end{definition}

We will see after stating Theorem~\ref{theo.main} that the functions
$h\mapsto\vol_{\CA^+}(h)$ and $h\mapsto N_{\CA^+}(h)$ are
respectively polynomial and quasipolynomial on each chamber of
$\CC(\CA)$.

The following formulae (see, for example,~\cite{ver}) compute  the
Laplace transform of the locally polynomial function
$\vol_{\CA^+}(h)$ and  the discrete Laplace transform of the
quasipolynomial function $N_{\CA^+}(h)$.

\begin{proposition} \label{prop.easy}
Let $u\in\CC(\CA)^*$. Then:
\begin{enumerate}
\item
  $\displaystyle \int_{\CC(\CA)}e^{-\ll h,u\rr}\vol_{\CA^+}(h)dh
    =\frac{1}{\prod_{i=1}^N\alpha_i(u)}$.
\item
  $\displaystyle \sum_{h\in V_{\Z}\cap\CC(\CA)}
      e^{-\ll h,u\rr}N_{\CA^+}(h)
  =\frac{1}{\prod_{i=1}^N(1-e^{-\ll\alpha_i,u\rr})}$.
\end{enumerate}
\end{proposition}

\section{Jeffrey-Kirwan residue}\label{sect.def.JK}

The aim of this section is  to explain some theoretical results
due to Jeffrey and Kirwan which are fundamental for our work. They
described an efficient scheme for computing the inverse Laplace
transforms in the context of hyperplane arrangements.

Let's go back to the space of rational functions $\CR_{\CA}$. It
is $\Z$-graded by degree. Of great importance for our exposition
will be certain functions in $\CR_{\CA}$ of degree $-r$. Every
function in $\CR_{\CA}$ of degree $-r$ may be decomposed into a
sum of basic fractions $f_\sigma$ (see Equation~\eqref{equa.basi})
and degenerate fractions; degenerate fractions are those for which
the linear forms in the denominator do not span $V$.
Given $\sigma\in\bases(\CA)$, we write $\CC(\sigma)$ for the
cone generated by $\alpha_i$ ($i\in\sigma$) and by
$\vol(\sigma)>0$ for the volume of the parallelotope
$\sum_{i=1}^r[0,1]\alpha_i$ computed for the measure $dh$.
Observe that $\vol(\sigma)=|\det(\sigma)|$, where $\sigma$ is the
matrix which columns are the $\alpha_i$'s.
Now having
fixed a chamber $\gc$, we define a functional $\JK_\gc(\phi)$ on
$\CR_{\CA}$ called the Jeffrey-Kirwan residue (or
\emph{JK residue}) as follows. Let
\begin{equation} \label{defi.jk}
\JK_{\gc}(f_\sigma)=
\begin{cases}
  \vol(\sigma)^{-1},&\mbox{if }\gc\subset\CC(\sigma),\\
  0, &\mbox{if }\gc\cap\CC(\sigma)=\emptyset.
\end{cases}
\end{equation}

By setting the value of the JK residue of a degenerate fraction or
that of a rational function of pure degree different from $-r$
equal to zero, we have defined the JK residue on $\CR_{\CA}$.

We may go further and extend the definition to the space
${\widehat\CR}_{\CA}$ which is the space consisting of functions
$P/Q$ where $Q$ is a product of powers of the linear forms
$\alpha_i$ and $P=\sum_{k=0}^{\infty}P_k$ is a formal power
series.
Indeed suppose that $P/Q\in{\widehat\CR}_{\CA}$ where we may
assume that $Q$ is of degree $q$, and
$P=\sum_{k=0}^{\infty}P_k$ is a formal power series with $P_k$ of
degree $k$.
Then we just define
$$\JK_{\gc}(P/Q)=\JK_{\gc}(P_{q-r}/Q)$$
as the JK residue of the component of degree $-r$ of $P/Q$.
In particular if $\phi\in\CR_{\CA}$ and $h\in V$, the function
$$e^{\ll h,u\rr}\phi(u)
  =\sum_{k=0}^{\infty}\frac{\ll h,u\rr^k}{k!}\phi(u)$$
is in ${\widehat\CR}_{\CA}$ and we may compute its JK
residue.
Observe that the JK residue depends on the measure $dh$.

Let's now make a short digression that should clarify why
JK residues compute inverse Laplace transforms.
For $u\in\CC(\CA)^*$ we have:
$$\frac{1}{\vol(\sigma)}\int_{\CC(\sigma)}
e^{-\ll h,u\rr}dh=f_\sigma(u).$$
In other words the inverse Laplace transform of
$f_\sigma$ computed at the point $h\in\CC(\CA)$ is
$\frac{1}{\vol(\sigma)}\chi_\sigma(h)$, where $\chi_\sigma$ is the
characteristic function of the cone $\CC(\sigma)$. We state this
as a formula:
$$\frac{1}{\vol(\sigma)}\chi_\sigma(h)=L^{-1}(f_{\sigma})(h).$$

Since the JK residue can be written in terms of basic fractions,
the following theorem~\cite{jeffreykirwan} is not surprising:

\begin{theorem}[Jeffrey-Kirwan] \label{theo.jk}
If $\phi\in\CR_{\CA}$, then for any $h\in\gc$ we have:
$$(L^{-1}\phi)(h)
  =\JK_{\gc}\left(e^{\ll h,\,\cdot\,\rr}\phi\right).$$
\end{theorem}

Assume that $\Psi:U\to U$ is a holomorphic transformation defined
on a neighborhood of $0$ in $U$ and invertible. We also assume
that $\alpha_j(F(u))=\alpha_j(u) f_j(u)$, where $f_j(u)$ is
holomorphic in a neighborhood of $0$ and $f_j(0)\neq 0$.

If $\phi$ is a function in ${\widehat\CR}_{\CA}$, the function
$\Psi^*\phi(u)=\phi(\Psi(u))$ is again in ${\widehat\CR}_{\CA}$.
Let $\jac(\Psi)$ be the Jacobian of the map $\Psi$.
The function $\jac(\Psi)$ is calculated as follows: write
$\Psi(u)
 =(\Psi_1(u_1,u_2,\ldots,u_r),\ldots,\Psi_r(u_1,u_2,\ldots,u_r))$.
Then
$\jac(\Psi)(u)
  =\det((\frac{\partial}
              {\partial u_i}\Psi_j)_{i,j})$.
We assume that $\jac(\Psi)(u)$ does not vanish at $u=0$. For any
$\phi$ in ${\widehat\CR}_{\CA}$ the following change of variable
formula \cite{baldonivergne}, Theorem 45, which will be useful in
our calculations later on, holds:

\begin{proposition} \label{prop.chng}
The Jeffrey-Kirwan residue obeys the rule of change of variables:
$$\JK_{\gc}(\phi)=\JK_{\gc}(\jac(\Psi)(\Psi^*\phi)).$$
\end{proposition}

We conclude this section by recalling  the formula for
$N_{\CA^+}(h)$. 

Consider the dual lattice
$U_\Z=\{u\in U\,|\,\ll u,V_\Z\rr\subset\Z\}$ and the torus
$T=U/U_\Z$.
Choosing a basis $\{u_1,\ldots,u_r\}$ of $U_\Z$ we may identify
$T$ with the subset of $U$ defined by the fundamental domain for
translation by $U_\Z$:
$$\left\{\sum_{j=1}^rt_j u_j\right\}$$
with $0\leq t_j<1$.

Every element $g$ in $T=U/U_\Z$ produces a function on $V_\Z$ by
$h\mapsto e^{\ll h,2\pi\sqrt{-1}G\rr}$, where we denote by $G$ a
representative of $g\in U/U_\Z$.\footnote{
We prefer to denote the complex number $i$ by $\sqrt{-1}$ because
we use $i$ for many indices.}
For $\sigma\in\bases(\CA)$ we
denote by $T(\sigma)$ the subset of $T$ defined by
$$T(\sigma)=\left\{g\in T\,\Big|\,\,
  e^{\ll\alpha,2\pi\sqrt{-1}G\rr}=1\,\,
  \mbox{for all}\,\alpha\in\sigma\right\}.$$
This is a finite subset of $T$. In particular if $\sigma$ is a
$\Z$-basis of $V_\Z$, then $T(\sigma)$ is reduced to the identity.
More generally, consider the lattice $\Z\sigma$ generated by the
elements $\alpha$ in $\sigma$. If $p$ is an integer such that
$\Z\sigma\subset p V_\Z$, then all elements of $T(\sigma)$ are of
order $p$.

For $g\in T$ and $h\in V_\Z$, consider the \emph{Kostant function}
$F(g,h)$ on $U$ defined by
\begin{equation}\label{equa.Fgh}
F(g,h)(u)=\frac{e^{\ll h,2\pi\sqrt{-1}G+u\rr}}
     {\prod_{i=1}^N(1-e^{-\ll\alpha_i,2\pi\sqrt{-1}G+u\rr})}.
\end{equation}
For example when $g=0$,
$$F(0,h)(u)=
\frac{e^{\ll h,u\rr}}
     {\prod_{i=1}^N(1-e^{-\ll\alpha_i,u\rr})}.$$
The function $F(g,h)(u)$ is an element of ${\widehat\CR}_{\CA}$.
Indeed if we write
$$I(g)=\left\{i\,\Big|\,\,1\leq i\leq N,
  e^{-\ll\alpha_i,2\pi\sqrt{-1}G\rr}=1\right\},$$
then
\begin{equation} \label{equa.fpsi}
F(g,h)(u)
  =e^{\ll h,2\pi\sqrt{-1}G\rr}\frac{e^{\ll h,u\rr}}
        {\prod_{i\in I(g)}\ll\alpha_i,u\rr}
   \psi^g(u)
\end{equation}
where $\psi^g(u)$ is the holomorphic function of $u$ (in a
neighborhood of zero) defined by
$$\psi^g(u)=
  \prod_{i\in I(g)}
    \frac{\ll\alpha_i,u\rr}
         {(1-e^{-\ll\alpha_i,u\rr})}
  \times
  \prod_{i\notin I(g)}
    \frac{1}
         {(1-e^{-\ll\alpha_i,2\pi\sqrt{-1}G+u\rr})}.$$
If $\gc$ is a chamber of $\CC(\CA)$, the Jeffrey-Kirwan
residue $\JK_\gc(F(g,h))$ is well defined.

The following theorem is due to Szenes-Vergne~\cite{SzeVer}. If
the set $\CA$ is unimodular (that is, each $\sigma \in
\bases(\CA)$ is a $\Z$-basis of $V_\Z$), it is a reformulation of
Khovanskii-Pukhlikhov Riemann-Roch calculus on simple
polytopes~\cite{KhoPuk}. For a general set $\CA$, this refines the
formula of Brion-Vergne \cite{BriVer97}.

\begin{theorem} \label{theo.main}
Let $\gc$ be a chamber of the cone $\CC(\CA)$ and
$\overline{\gc}$ its closure.
Then:
\begin{enumerate}
\item For $h\in\overline{\gc}$ we have
  $$\vol_{\Z,\CA^+}(h)
    =\vol(V/V_\Z,dh) \JK_{\gc}\left(\frac{e^{\ll h,\,\cdot\,\rr}}
                        {\prod_{i=1}^N\alpha_i}\right).$$
\item\label{theo.main.two} Assume that $F$ is a finite subset of $T$ such that for any
  $\sigma\in\bases(\CA)$, we have $T(\sigma)\subset F$.
  Then for $h\in V_\Z\cap\overline{\gc}$, we have
  $$N_{\CA^+}(h)=\vol(V/V_\Z,dh)\sum_{g\in F}\JK_{\gc}(F(g,h)).$$
\end{enumerate}
\end{theorem}

Observe that the right-hand side of (\ref{theo.main.two}) does
not depend on the measure $dh$, as it should be.

Let us explain the behavior of these functions on a chamber $\gc$.
By definition, a quasipolynomial function on a lattice $L$ is a
linear combination of products of polynomial functions and of
periodic functions (functions constants on cosets $h+p L$ where
$p$ is an integer). We now show that the normalized volume
$\vol_{\Z,\CA^+}(h)$ is given by a polynomial formula, when $h$
varies in a chamber $\overline \c$, while $N_{\CA^+}(h)$ is given
by a quasipolynomial formula when $h$ varies in
$V_\Z\cap\overline{\gc}$.

The residue vanishes except on degree $-r$, so that
$$\JK_{\gc}\left(\frac{e^{\ll h,u\rr}}
                        {\prod_{i=1}^N\ll\alpha_i,u\rr}\right)
=\frac{1}{(N-r)!}\JK_{\gc}\left(\frac{\ll h,u\rr^{N-r}}
                        {\prod_{i=1}^N\ll\alpha_i,u\rr}\right) , $$
and as expected the normalized volume is a polynomial
homogeneous function of $h$ of degree $N-r$ on each chamber.

Now if $\CA$ is unimodular then the above defined set $F$
is $\{0\}$.
Hence the number of integral points in the polytope
$\Pi_{\CA^+}(h)$ satisfies
$$N_{\CA^+}(h)= \vol \left( V/V_\Z, dh \right)
\JK_{\gc}\left(\frac{e^{\ll h,u\rr}}
{\prod_{i=1}^N(1-e^{-\ll\alpha_i,u\rr})}\right)$$
and is a polynomial of degree $N-r$ whose homogeneous
component of degree $N-r$ is the normalized volume. More
precisely, write
$$\frac{e^{\ll h,u\rr}}
       {\prod_{i=1}^N(1-e^{-\ll\alpha_i,u\rr})}
= \frac{e^{\ll h,u\rr}}
       {\prod_{i=1}^N\ll\alpha_i,u\rr}
  \times
  \frac{\prod_{i=1}^N\ll\alpha_i,u\rr}
       {\prod_{i=1}^N(1-e^{-\ll\alpha_i,u\rr})}$$
where
$$\frac{\prod_{i=1}^N\ll\alpha_i,u\rr}
       {\prod_{i=1}^N(1-e^{-\ll\alpha_i,u\rr})}
= \sum_{k=0}^{+\infty}\psi_k(u)$$
is a holomorphic function of $u$ in a neighborhood of $0$ with
$\psi_0(u)=1$. Consequently
\begin{eqnarray}
N_{\CA^+}(h)
&=&\vol \left( V/V_\Z, dh \right) \JK_{\gc}\left(
           \frac{e^{\ll h,u\rr}}
                {\prod_{i=1}^N\ll\alpha_i,u\rr}
           \times
           \sum_{k=0}^{+\infty}\psi_k(u)
           \right) \nonumber \\
&=&\vol \left( V/V_\Z, dh \right) \sum_{k=0}^{N-r}
   \frac{1}{(N-r-k)!}
   \JK_{\gc}\left(
           \frac{\ll h,u\rr^{N-r-k}\psi_k(u)}
                {\prod_{i=1}^N\ll\alpha_i,u\rr}
           \right). \label{na+formula3}
\end{eqnarray}
The unimodular case applies to the root  system $A_r$.

Finally in the non-unimodular case (for example for the root
systems $B_r$, $C_r$, $D_r$) the set $F$ is no longer reduced to
$\{0\}$. Let us denote by
$\psi^g(u)=\sum_{k=0}^{+\infty}\psi_k^g(u)$ the series development
of the holomorphic function $\psi^g$ appearing in
formula~\eqref{equa.fpsi}. Then we see that $\JK_{\gc}(F(g,h))$
equals
\begin{eqnarray}
& &\JK_{\gc}
       \left(e^{\ll h,2\pi\sqrt{-1}G\rr}\frac{e^{\ll h,u\rr}}
                  {\prod_{i\in I(g)}\ll\alpha_i,u\rr}
             \psi^g(u)
       \right) \label{na+formula1} \\
&=&e^{\ll h,2\pi\sqrt{-1}G\rr}\sum_{k=0}^{|I(g)|-r}
      \frac{1}{(|I(g)|-r-k)!}
      \JK_{\gc}\left(\frac{\ll h,u\rr^{|I(g)|-r-k}}
                         {\prod_{i\in I(g)}\ll\alpha_i,u\rr}
                    \psi_k^g(u)
              \right). \nonumber
\end{eqnarray}
If $g$ is of order $p$, the function $h\mapsto e^{\ll
h,2\pi\sqrt{-1}G\rr}$ is constant on each coset $h+pV_\Z$ of the
lattice $pV_\Z$, while the function $h\mapsto
\JK_{\gc}\left(\frac{\ll h,u\rr^{|I(g)|-r-k}}
                         {\prod_{i\in I(g)}\ll\alpha_i,u\rr}
                    \psi_k^g(u)
              \right)$
 is a polynomial function of $h$ of degree $|I(g)|-r-k$.
 Thus the function
\begin{equation}\label{na+formula2}
 N_{\CA^+}(h)=\vol \left( V/V_\Z, dh \right) \sum_{g\in F}\JK_{\gc}(F(g,h))
\end{equation}
 is given by a quasipolynomial formula when $h$ varies in the closure
 of a chamber. Note that its highest degree component is polynomial and is the
normalized volume as expected.

\begin{example} \label{exam.form.B2}
Let us compute the normalized volume and number of integral points
for the root system $B_2$, that is for
$\CA^+=\CB_2=\{e_1,e_2,e_1+e_2,e_1-e_2\}$. Fix a chamber $\gc$ and
an integral vector $h=(h_1,h_2)$ in the cone $\CC(\CB_2)$.
Observe that the root lattice is $\Z e_1 \oplus \Z e_2 $ and
$\vol \left( V/V_\Z, dh \right) = 1$ for the measure $dh=dh_1 dh_2$.
Then the normalized volume equals
\begin{eqnarray*}
&&\frac{1}{2!}
   \JK_{\gc}\left(\frac{(h_1u_1+h_2u_2)^2}{u_1u_2(u_1+u_2)(u_1-u_2)}
           \right).
\end{eqnarray*}
Note that
\begin{eqnarray*}
\frac{u_1^2}{u_1u_2(u_1+u_2)(u_1-u_2)}
 &=&\frac{1}{u_2(u_1+u_2)}+\frac{1}{(u_1+u_2)(u_1-u_2)}
\end{eqnarray*}
(and similar quotients of $u_2^2$ and $u_1u_2$),
so that the normalized volume is
\begin{eqnarray*}
\frac12 \JK_{\gc}\left(
                 \frac{h_1^2}{u_2(u_1+u_2)}
                +\frac{h_1^2+2h_1h_2+h_2^2}{(u_1+u_2)(u_1-u_2)}
                -\frac{h_2^2}{u_1(u_1+u_2)}
                \right).
\end{eqnarray*}
There are three chambers, namely
$\gc_1=\CC(\{e_2,e_1+e_2\})$, $\gc_2=\CC(\{e_1,e_1+e_2\})$,
$\gc_3=\CC(\{e_1-e_2,e_1\})$ (see Figure~\ref{figu.B2}). Now let us
compute the Jeffrey-Kirwan residues on the chambers. As
$$\begin{array}{rclrcl}
\JK_{\gc_1}\left(\frac{1}{u_2(u_1+u_2)}\right)&=&1,&
\JK_{\gc_2}\left(\frac{1}{(u_1+u_2)(u_1-u_2)}\right)&=&\frac12,\\
\JK_{\gc_2}\left(\frac{1}{u_1(u_1-u_2)}\right)&=&1,&
\JK_{\gc_3}\left(\frac{1}{(u_1+u_2)(u_1-u_2)}\right)&=&\frac12,
\end{array}$$
we obtain
\begin{eqnarray*}
\vol(\Pi_{\CB_2}(h))&=&\frac12 h_1^2\quad\mbox{ if }h\in\gc_1,\\
\vol(\Pi_{\CB_2}(h))&=&\frac14 (h_1+h_2)^2-\frac 12 h_2^2\quad\mbox{ if }h\in\gc_2,\\
\vol(\Pi_{\CB_2}(h))&=&\frac14 (h_1+h_2)^2\quad\mbox{ if }h\in\gc_3.\\
\end{eqnarray*}
Note that the formulae agree on walls $\gc_1\cap\gc_2$
and $\gc_2\cap\gc_3$.

\begin{figure}[t]
  \vbox{\epsfxsize=4cm
    \epsfbox{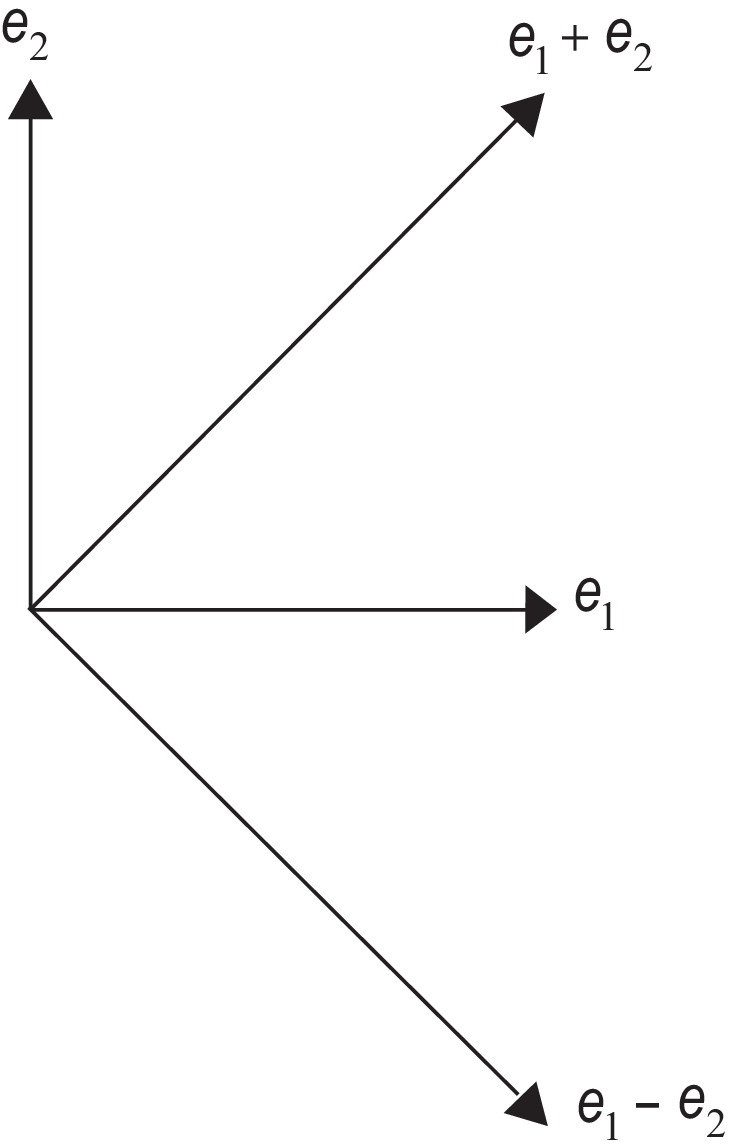}
    \vskip -5.5cm \hskip -1.3cm {\large$\mathfrak{c}_1$}
    \vskip 0.5cm \hskip 0.8cm {\large$\mathfrak{c}_2$}
    \vskip 1.8cm \hskip 0.8cm {\large$\mathfrak{c}_3$}
  }
\vspace{2cm}
\caption{The 3 chambers for $B_2$}
\label{figu.B2}
\end{figure}


%
For the number of integral points, we first note that
$F=\{(0,0),(1/2,1/2)\}$.
Consequently $N_{\CB_2}(h)$ is equal to the Jeffrey-Kirwan residue
of $f_1=F((0,0),h)$ plus $f_2=F((1/2,1/2),h)$.
We rewrite the series $f_j$ ($j=1$, $2$) as
$f_j=f'_j\times e^{u_1h_1+u_2h_2}/u_1u_2(u_1+u_2)(u_1-u_2)$
where
\begin{eqnarray*}
f'_1&=&\frac{u_1}{1-e^{-u_1}}
       \times
       \frac{u_2}{1-e^{-u_2}}
       \times
       \frac{u_1+u_2}{1-e^{-(u_1+u_2)}}
       \times
       \frac{u_1-u_2}{1-e^{-(u_1-u_2)}},\\
f'_2&=&\frac{u_1}{1+e^{-u_1}}
       \times
       \frac{u_2}{1+e^{-u_2}}
       \times
       \frac{u_1+u_2}{1-e^{-(u_1+u_2)}}
       \times
       \frac{u_1-u_2}{1-e^{-(u_1-u_2)}}\times(-1)^{h_1+h_2}.
\end{eqnarray*}
Using the series expansions
$\frac{x}{1-e^{-x}}=1+\frac12x+\frac{1}{12}x^2+O(x^3)$ and
$\frac{x}{1+e^{-x}}=\frac12x+O(x^2)$, we obtain that the number of
integral points is the JK residue of
\begin{eqnarray*}
&& \frac{u_1(1+\frac32h_1+\frac12h_1^2)}{u_2(u_1-u_2)(u_1+u_2)}
  +\frac{\frac34+h_1h_2+\frac32h_2+\frac12h_1}{(u_1-u_2)(u_1+u_2)}
  +\frac{u_2(\frac12h_2^2+\frac12h_2)}{u_1(u_1-u_2)(u_1+u_2)}\\
&&+(-1)^{h_1+h_2}\frac{\frac14}{(u_1+u_2)(u_1-u_2)}\\
\end{eqnarray*}
\begin{eqnarray*}
&& =\frac{(1+\frac32h_1+\frac12h_1^2)}{u_2(u_1+u_2)}
  +\frac{\frac74+2(h_1+h_2)+\frac12h_2^2+h_1h_2+\frac12h_1^2}{(u_1-u_2)(u_1+u_2)}
  -\frac{(\frac12h_2^2+\frac12h_2)}{u_1(u_1+u_2)}\\
&&+(-1)^{h_1+h_2}\frac{\frac14}{(u_1+u_2)(u_1-u_2)}.\\
\end{eqnarray*}
We then obtain:
\begin{eqnarray*}
N_{\CB_2}(h)&=& 1+\frac32h_1+\frac12h_1^2\\
             & & \quad\mbox{ if }h\in\gc_1,\\
N_{\CB_2}(h)&=&
\frac14h_1^2+\frac12h_1h_2-\frac14h_2^2+h_1+\frac12h_2+\frac78+(-1)^{h_1+h_2}\frac18\\
             & &
              \quad\mbox{ if }h\in\gc_2,\\
N_{\CB_2}(h)&=&\frac14h_1^2+\frac12h_1h_2+\frac14h_2^2+h_1+h_2+\frac78+(-1)^{h_1+h_2}\frac18\\
             & & \quad\mbox{ if }h\in\gc_3.
\end{eqnarray*}
Note that the functions $N_{\CB_2}$ agree on walls, and the
formulae above are valid on the closures of the chambers.

Our general method to implement Theorem~\ref{theo.main}
for root systems is more systematic and will be explained in the
course of this article.
\end{example}

\begin{remark}
Combining \eqref{na+formula1} and \eqref{na+formula2}, we can see
that the quasipolynomial character of the integral-point counting
functions $N_\CA^+$ stems precisely from the root of unity in
\eqref{na+formula1}. Furthermore, we will see in Lemmas
\ref{lemm.B.F}, \ref{lemm.C}, and \ref{lemm.D.F} that for root
systems of type $B$, $C$, and $D$, these roots of unity are of
order 2, as in the above example for $\CB_2$. (For root systems of
type $A$, \eqref{na+formula3} shows that $N_\CA^+$ is always a
polynomial.) Let us record the following immediate consequence:
\end{remark}

\begin{corollary}
The integral-point counting functions $N_{ \CB_r } , N_{ \CC_r } ,
N_{ \CD_r }$ are quasipolynomials with period $2$.
\end{corollary}

\begin{remark}
The partition functions $N_{ \CA_r } , N_{ \CB_r } , N_{ \CC_r } ,
N_{ \CD_r }$ can be interpreted as (weak) flow quasipolynomials on
certain signed graphs \cite{BeckZas}. The polynomiality of $N_{
\CA_r }$ follows immediately from this interpretation and a
unimodularity argument; the fact that the quasipolynomials $N_{
\CB_r } , N_{ \CC_r } , N_{ \CD_r }$ have period 2 follows from a
half-integrality result of Lee \cite{Lee}.
\end{remark}

\begin{remark} \label{rema.unim}
In the case where $\CA^+$ is an arbitrary sequence of vectors in
$V_\Z$, the straightforward implementation of  Theorem
\ref{theo.main} above is of exponential complexity. Indeed we make
a summation on the set $F$, which can become arbitrarily large.
Barvinok uses a signed cone decomposition to obtain an algorithm
of polynomial complexity, when the number of elements of $\CA^+$
is fixed, to compute the number $N_{\CA^+}(h)$; the {\tt LattE}
team implemented Barvinok's algorithm~\cite{lattemanual,latte} in
the language C. Our work will be dealing either with volumes of
polytopes, where the set $F$ does not enter, or with partition
function of classical root systems, where the set $F$ is
reasonably small. Then we obtain a fast algorithm, implemented for
the moment in the formal calculation software {\sc Maple}. This
algorithm for these particular cases can reach examples not
obtainable by the {\tt LattE} program.
\end{remark}

In next Section~\ref{sect.JK} we will give the basic formula for
$\JK_{\gc}$ involving maximal proper nested sets, as developed in
\cite{deconciniprocesi}, and iterated residues. These formulae are
implemented in our algorithms.

\section{A formula for the Jeffrey-Kirwan residue} \label{sect.JK}

If $f$ is a meromorphic function of one variable $z$ with a pole
of order less than or equal to $h$ at $z=0$ then we can write
$f(z)=Q(z)/z^h$, where $Q(z)$ is a holomorphic function near $z=0$.
If the Taylor series of $Q$ is given by
$Q(z)=\sum_{k=0}^{\infty}q_k z^k$, then as usual the residue at
$z=0$ of the function $f(z)=\sum_{k=0}^{\infty}q_k z^{k-h}$ is
the coefficient of $1/z$, that is, $q_{h-1}$.
We will denote it by $\res_{z=0}f(z)$.
To compute this residue we can either expand $Q$ into a power series
and search for the coefficient of $z^{-1}$, or employ the formula
\begin{equation} \label{equa.resf}
\res_{z=0}f(z)
=\frac{1}{(h-1)!}(\partial_z)^{h-1}\left(z^hf(z)\right)\Big|_{z=0}.
\end{equation}

We now introduce the notion of iterated residue on the space
$\CR_{\CA}$.

Let $\vec{\nu}=[\alpha_1,\alpha_2,\ldots,\alpha_r]$ be an ordered
basis of $V$ consisting of elements of $\CA$ (here we have
implicitly renumbered the elements of $\CA$ in order that the
elements of our basis are listed first). We choose a system of
coordinates on $U$ such that $\alpha_i(u)=u_i$. A function
$\phi\in\CR_{\CA}$ is thus written as a rational fraction
$\phi(u_1,u_2,\ldots,u_r)
  =\frac{P(u_1,u_2,\ldots,u_r)}{Q(u_1,u_2,\ldots,u_r)}$
where the denominator $Q$ is a product of linear forms.

\begin{definition} \label{defi.ires}
If $\phi\in\CR_{\CA}$, the iterated residue
$\Ires_{\vec{\nu}}(\phi)$ of $\phi$ for $\vec\nu$ is the scalar
$$\Ires_{\vec{\nu}}(\phi)
  =\res_{u_r=0}\res_{u_{r-1}=0}
   \cdots\res_{u_1=0}\phi(u_1,u_2,\ldots,u_r)$$
where each residue is taken assuming that the variables
with higher indices are considered constants.
\end{definition}

Keep in mind that at each step the residue operation augments the
homogeneous degree of a rational function  by $+1$ (as for example
$\res_{x=0}(1/xy)=1/y$) so that the iterated residue vanishes on
homogeneous elements $\phi\in \CR_{\CA}$, if the homogeneous
degree of $\phi$ is different from $-r$.

Observe that the value of $\Ires_{\vec{\nu}}(\phi)$ depends on the
order of ${\vec{\nu}}$. For example, for $f=1/(x(y-x))$ we have
$\res_{x=0}\res_{y=0}(f)=0$ and $\res_{y=0}\res_{x=0}(f)=1$.

\begin{remark} \label{rema.wedg}
Choose any basis $\gamma_1$, $\gamma_2$, \ldots, $\gamma_r$ of $V$
such that $\oplus_{k=1}^j\alpha_j=\oplus_{k=1}^j\gamma_j$ for
every $1\leq j\leq r$  and such that
$\gamma_1\wedge\gamma_2\wedge\cdots\wedge\gamma_r
  =\alpha_1\wedge\alpha_2\wedge\cdots\wedge\alpha_r$.
Then, by induction, it is easy to see that for $\phi\in \CR_{\CA}$

$$\res_{\alpha_r=0}\cdots\res_{\alpha_1=0}\phi
=\res_{\gamma_r=0}\cdots\res_{\gamma_1=0}\phi.$$

Thus given an ordered basis, we may modify $\alpha_2$ by
$\alpha_2+c\alpha_1$, \ldots, with the purpose of getting easier
computations.
\end{remark}

The following lemma will be useful later on.

\begin{lemma}\label{necessary}
Let $\vec{\nu}=[\alpha_1,\alpha_2,\ldots,\alpha_r]$ and
$f_\beta=\frac{1}{\prod_{i=1}^r\beta_i}$ be a basic fraction. Then
 the  iterated residue $\Ires_{\vec{\nu}}(f_\beta)$ is non zero if
and only if there exists a permutation $w$ of $\{1,2,\ldots,r\}$
such that:
\begin{eqnarray*}
\beta_{w(1)}&\in&\R\alpha_1,\\
\beta_{w(2)}&\in&\R\alpha_1\oplus\R\alpha_2,\\
&\vdots&\\
\beta_{w(r)}
  &\in&\R\alpha_1\oplus\cdots\oplus\R\alpha_r.
\end{eqnarray*}

\end{lemma}

\begin{definition}\label{zefepsilon}
Let $\vec{\nu}=[\alpha_1,\alpha_2,\ldots,\alpha_r]$ and let
$u_j=\alpha_j(u)$. Choose a sequence of real numbers:
$0<\epsilon_1<\epsilon_2<\cdots<\epsilon_r$.
Then define the torus
\begin{equation}
  T({\vec{\nu}})
  =\{u\in U_\C\,|\,|u_j|=\epsilon_j,\,j=1,\ldots,r\}.
\end{equation}
\end{definition}
The torus $T({\vec{\nu}})$ is identified via the basis $\alpha_j$
with the product of $r$ circles oriented counterclockwise. The sequence
$\left[\epsilon_1, \epsilon_2, \dots, \epsilon_r\right]$
is chosen so that elements $\alpha_q$ not in
$\oplus_{k=1}^j\R\alpha_j$ do not vanish on the domain
$\{u\in U_\C\,|\,|u_k|\leq\epsilon_k,1\leq k\leq j\,;\,
                 |u_i|=\epsilon_i,i=j+1,\ldots,r\}$.
This is achieved by choosing the ratios
$\epsilon_j/\epsilon_{j+1}$ very small. The torus $T({\vec{\nu}})$
is contained in $U_{\C}(\CA)$ and the homology class
$[T({\vec{\nu}})]$ of this torus is independent of the choice of
the sequence of the ordered $\epsilon_j$ \cite{SzeVerToric}.

Choose  an ordered basis $e_1,e_2,\ldots, e_r$ of $V$ of volume
$1$ with respect to the measure $dh$. For $z\in U_\C$, define
$z_j= \ll z, e_j \rr $ and  $dz=dz_1\wedge dz_2\wedge\cdots \wedge
dz_r$. Denote by $\det({\vec{\nu}})$ the determinant of the basis
$\alpha_1,\alpha_2,\ldots,\alpha_r$ with respect to  the basis
$e_1,e_2,\ldots,e_r$.

\begin{lemma}\label{lemm.inte}
For $\phi\in\CR_{\CA}$, we have
$$\frac{1}{\det({\vec{\nu}})}
  \res_{\alpha_r=0}\cdots\res_{\alpha_1=0}\phi
=\frac{1}{(2\pi\sqrt{-1})^r}\int_{T({\vec{\nu}})}\phi(z)dz.$$
\end{lemma}

Thus, as for the usual residue, the iterated residue  can be
expressed as an integral.


We now introduce the notion of maximal proper nested set,
MPNS in short.

De Concini-Procesi~\cite{deconciniprocesi} prove that the set of
MPNS is in bijection with  the so-called no broken circuits bases
of $\CA$ (with respect to a order to be specified).
This is helpful as the JK residue can be computed in terms of
iterated residues with respect to these bases.

If $S$ is a subset of $\CA$, we denote by $\ll S\rr$ the vector
space spanned by $S$. More generally if $M=\{ S_i \}$ is a set of
subsets of $\CA$, we denote by $\ll M \rr$ the vector space
spanned by all elements of the sets $S_i$. We say that a subset
$S$ of $\CA$ is \emph{complete} if $S=\ll S\rr\cap\CA$ or in other
words if any linear combination of elements of $S$ belongs to $S$.
A complete subset $S$ is called \emph{reducible} if we can find a
decomposition $V=V_1\oplus V_2$ such that $S=S_1\cup S_2$ with
$S_1\subset V_1$ and $S_2\subset V_2$. Otherwise $S$ is said to be
\emph{irreducible}.

\begin{definition} \label{defi.ns}
Let $\CI$ be the set of irreducible subsets of $\CA$. A set
$M=\{I_1,I_2,\ldots,I_k\}$ of irreducible subsets of $\CA$ is
called \emph{nested} if, given any subfamily $\{I_1,\ldots,I_m\}$ of $M$
such that there exists no $i$, $j$ with $I_i\subset I_j$, then the
set $I_1\cup\cdots\cup I_m$ is {\bf complete} and the elements
$I_j$ are the irreducible components of $I_1\cup I_2\cup\cdots\cup
I_m$.
\end{definition}

\begin{example} \label{exam.Kn}
Let $E$ be an $r+1$-dimensional vector space with basis $e_i$
($i=1$, \ldots, $r$). We consider the set
$$\CK_{r+1}=\{e_i-e_j\,|\,1\leq i<j\leq r+1\}.$$
These are the positive roots for the system $A_{r}$. The
irreducible subsets of $\CK_{r+1}$ are indexed by subsets $S$ of
$\{1,2,\ldots,r+1\}$, the corresponding irreducible subset being
$\{e_i-e_j\,|\,i,j\in S, i<j\}$. For instance the set
$S=\{1,2,4\}$ parametrizes the set of roots given by
$\{e_1-e_2,e_2-e_4,e_1-e_4\}$.

A nested set is represented by a collection
$M=\{S_1,S_2,\ldots,S_k\}$ of subsets of $\{1,2,\ldots,r+1\}$ such
that if $S_i$, $S_j\in M$ then either $S_i\cap S_j$ is empty, or
one of them is contained in another.
\end{example}

\begin{definition} \label{defi.mns}
A maximal nested set (in short MNS) $M$ is a nested set such that
for every irreducible set $\CI$ of $\CA$ the set $M\cup \{\CI\}$
is no longer nested.
\end{definition}

A maximal nested set has exactly $r$ elements~\cite{deconciniprocesi}.

 Assume now that
$\CA$ is irreducible, otherwise just take the irreducible
components. Then every maximal nested set $M$ contains $\CA$. Let
$I_1$, $I_2$, \ldots, $I_k$ be the maximal elements of the set
$M\setminus\CA$. We see that the vector space spanned by $\ll
I_1\rr\oplus\ll I_2\rr\oplus\cdots\oplus\ll I_k\rr$ is of
codimension $1$ \cite[Proposition 1.3]{deconciniprocesi}.

\begin{definition} \label{defi.H.adm}
A hyperplane $H$ in $V$ is \emph{$\CA$-admissible} if it
is spanned by a set of vectors of $\CA$.
\end{definition}

Thus if $M$ is a MNS, the vector space  $\ll M\setminus \CA\rr $
is an admissible hyperplane $H$.

\begin{definition} \label{defi.atta}
Let $\CA$ be irreducible  and let $H$ be a  $\CA$-admissible
hyperplane. All MNPS's such that $\ll M\setminus \CA\rr =H$  are
said attached to $H$.
\end{definition}

Therefore to classify maximal nested sets (MNS) for an irreducible
set  $\CA$ we proceed by running over the set of $\CA$-admissible
hyperplanes, as described in Figure~\ref{figu.algo.MNS}.

\begin{figure}[ht]
\begin{itemize}
\item Take a hyperplane $H$ spanned by a set of vectors of
  $\CA$.
\item Break $\CA\cap H$ into irreducible subsets
  $I_1\cup I_2\cup\cdots\cup I_k$.
\item For each irreducible $I_i$ construct the set
  $\{M^i_{1},\ldots,M^i_{k_i}\}$ of maximal nested sets for
  $I_i$.
\item Set $C_i=\{1,\ldots,k_i\}$.
\item A maximal nested set is then given by the union
  $M^1_{c_1}\cup M^2_{c_2}\cup\cdots\cup M^k_{c_k}\cup\{\CA\}$
   where $c_1\in C_1$, \ldots, $c_k\in C_k$,
  and all of them are obtained by letting $c_i$ vary.
\end{itemize}
\caption{Building of all MNSs attached to an $\CA$-admissible
hyperplane $H$} \label{figu.algo.MNS}
\end{figure}

The whole algorithm will be described in detail in
Figure~\ref{figu.algo.MPNS}, Section~\ref{sect.gene}.


We describe now the notion of maximal proper nested set of $\CA$.

Fix a total order on the set $\CA$.
For example, we can choose a linear functional $\h$ on $V$
so that the values $\h(\alpha_i)$ are all distinct and positive.
Thus the value $\h(\alpha)$ is larger if $\alpha$ is deeper in
the interior of the cone.

Let $M=\{S_1,S_2,\ldots,S_k\}$ be a set of subsets of $\CA$.
In each $S_j$ we choose the element $\alpha_j$ maximal for the
order given by $\h$.
This defines a map $\theta$ from $M$ to $\CA$.

\begin{definition} \label{defi.mpns}
A maximal nested set $M$ is called \emph{proper} if $\theta(M)$ is a
basis of $V$.
We denote by $\CP(\CA)$ the set of maximal proper nested sets,
in short MPNS.
\end{definition}

If $M=\{I_1,I_2,\ldots,I_r\}$ is a maximal nested set, we
associate to  $M$ the list
$[\theta(I_{i_1}),\ldots,\theta(I_{i_r})]$ using the total order
on the elements $\theta(M)$; that is we have
$\h(\theta(I_{i_1}))<\h(\theta(I_{i_2}))<\cdots<\h(\theta(I_{i_r}))$.
Observe that, if $\CA$ is irreducible,  for every maximal nested
set, $I_{i_r}$ is always equal to $\CA$ and $\theta(I_{i_r})$ is
the highest element of $\CA$. We  will often implicitly
renumber our elements in $M$ such that
$\h(\theta(I_{1}))<\h(\theta(I_{2}))<\cdots<\h(\theta(I_{r}))$.

So we have associated to every maximal proper nested set $M$ an
ordered basis
$\overrightarrow{\theta(M)}=[\alpha_1,\alpha_2,\ldots,\alpha_r]$
of elements of $\CA$. In all implementations, we calculate
$\overrightarrow{\theta(M)}$ from a MNS $M$ with the procedure
${\tt ThetaMNS(M)}$. We denote by $\vol(M)>0$ the volume of the
parallelepiped $\sum_{i=1}^r[0,1]\alpha_i$ with respect to our
measure, and by $\CC(M)=\sum_{i=1}^r\R_{\geq 0}
\alpha_i\subset\CC(\CA)$ the cone generated by $\theta(M)$.

If $v$ is a regular element of $V$, let
\label{equa.goodforv}
\begin{equation}
\CP(v,\CA)=\{M\in\CP(\CA)\,|\,v\in\CC(M)\}.
\end{equation}
The set $\CP(v,\CA)$ depends only of the chamber $\gc$ where $v$
belongs.
We are now ready to state the basic formula for our calculations.

\begin{theorem}[DeConcini-Procesi, \cite{deconciniprocesi}] \label{theo.DCP}
Let $\gc$ be a chamber and let $v\in\gc$. Then, for $\phi\in
\CR_{\CA}$, we have
$$\JK_\gc(\phi)
  =\sum_{M\in\CP(v,\CA)}
    \frac{1}
         {\vol(M)}\Ires_{\overrightarrow{\theta(M)}}\phi.$$
\end{theorem}

We will use also the corresponding integration formula.

Each maximal proper nested set $M\in\CP(v,\CA)$ determines an
oriented cycle
$\left[T\left(\overrightarrow{\theta(M)}\right)\right]$ contained
in the open set $U_{\C}(\CA)$, as described in
Definition~\ref{zefepsilon}.

\begin{definition} \label{defi.cycl}
Let $\gc$ be a chamber.
 Define the oriented cycle:
$$H(\c)=\sum_{M\in\CP(v,\CA)}\operatorname{sign}
  \left(\det\left(\overrightarrow{\theta(M)}\right)\right)
  \left[T\left(\overrightarrow{\theta(M)}\right)\right].$$
\end{definition}

The following integral version of Theorem \ref{theo.DCP} will be
useful.

\begin{theorem}\label{theo.intJK}
Let $\gc$ be a chamber. Then for $\phi\in \CR_{\CA}$ we have
$$\JK_\gc(\phi)
  =\frac{1}{(2\pi\sqrt{-1})^r}\int_{H(\c)}\phi(z) dz.$$
\end{theorem}

The following example should help clarifying the notions
introduced.

\begin{example}\label{exam.K4}
We consider the set $\CK_4$ of positive roots for $A_3$
(see Figure~\ref{figu.A3}) defined by
$$\CK_4=\{e_i-e_j\,|\,1\leq i<j\leq 4\}.$$

\begin{figure}[ht]
  \centering
  \psfrag{a1}{$e_1$}
  \psfrag{a2}{$e_2$}
  \psfrag{a3}{$e_3$}
  \psfrag{a12}{$e_1-e_2$}
  \psfrag{a23}{$e_2-e_3$}
  \psfrag{a13}{$e_1-e_3$}
  \psfrag{a.1}{$a_1\geq 0$}
  \psfrag{a.2}{$a_2\geq 0$}
  \psfrag{a.3}{$a_3\geq 0$}
  \psfrag{a.12}{$a_1+a_2\geq 0$}
  \psfrag{a.13}{$a_1+a_3\geq 0$}
  \psfrag{a.23}{$a_2+a_3\geq 0$}
  \psfrag{a.123}{$a_1+a_2+a_3\geq 0$}
\includegraphics[height=8cm,width=8cm]{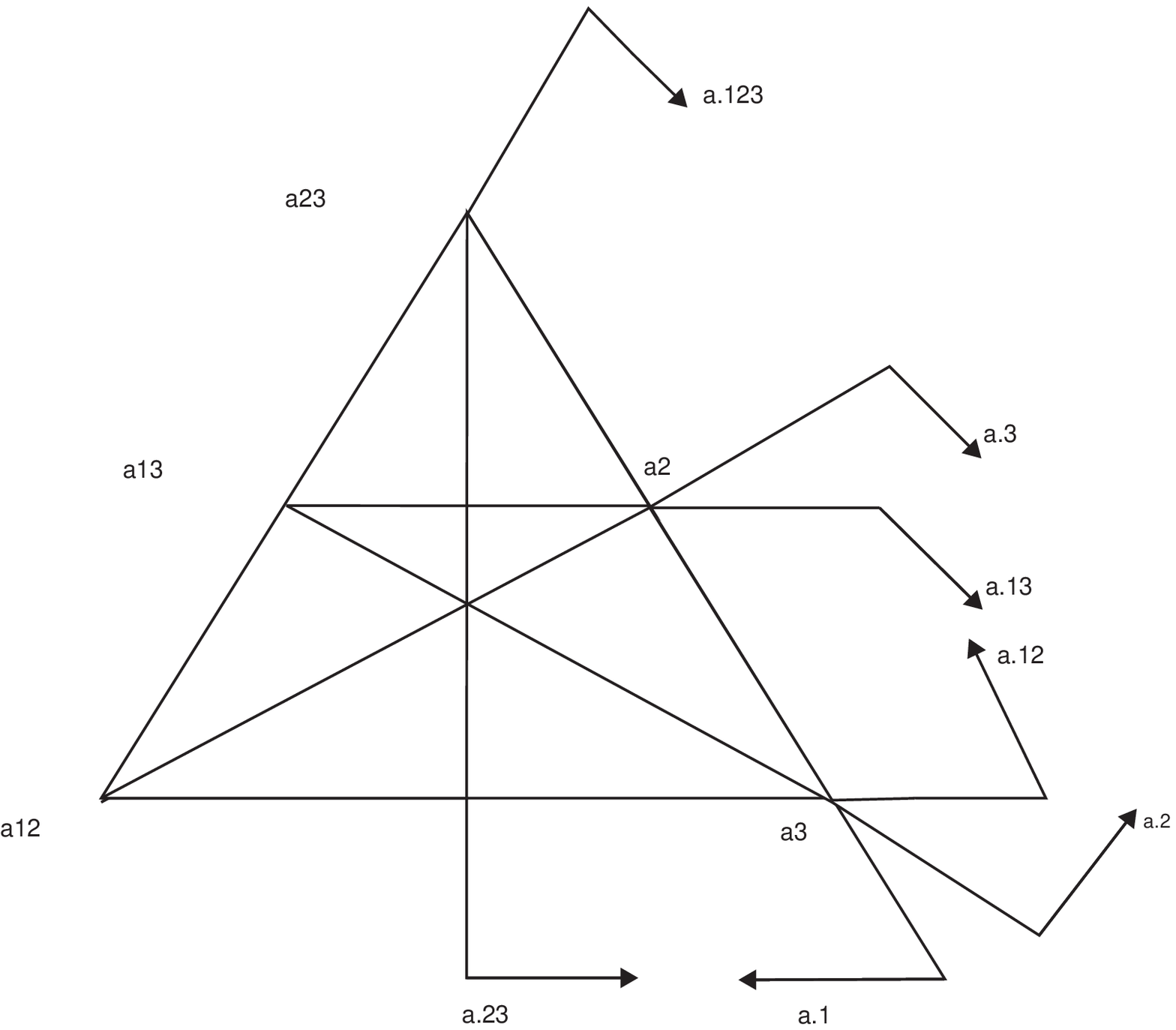}
\caption{Hyperplanes for $A_3$ with
$a=a_1e_1+a_2e_2+a_3e_3-(a_1+a_2+a_3)e_4$} \label{figu.A3}
\end{figure}
We let $V$ be the vector space generated by the elements in
$\CK_4$.
Then $V$ has dimension $3$  and we write an element of $V$ as
$$a=a_1 e_1+a_2 e_2+a_3 e_3-(a_1+a_2+a_3)e_4.$$
We consider the height function defined by
$$\h(e_1-e_2)=10,\quad\h(e_2-e_3)=11,\quad\h(e_3-e_4)=12.$$
This choice gives the following order on the roots:
$$e_1-e_2<e_2-e_3<e_3-e_4<e_1-e_3<e_2-e_4<e_1-e_4.$$
Take a hyperplane $H$ in $V$ spanned  by  two linearly
independent elements of $\CK_4$. Therefore it is the kernel of a
linear form $\sum_{i\in I_H}a_i$, where $I_H$ is a proper subset
of  $\{1,2,3,4\}$. The set of complementary indices gives the same
hyperplane. Thus each admissible  hyperplane partitions the set of
indices $\{1,2,3,4\}$ in two sets $Z_1$ and $Z_2$, where
$Z_1:=\{i\in I_H\}$ and $Z_2$ is the set of complementary indices.
In our example we have $7$ choices of admissible hyperplanes
corresponding to the following partitions:
$$\begin{array}{lll}
H_1=\{[1,2,3],[4]\},&H_2=\{[1,2,4],[3]\},&H_3=\{[1,3,4],[2]\},\\
H_4=\{[2,3,4],[1]\},&H_5=\{[1,2],[3,4]\},&H_6=\{[1,3],[2,4]\},\\
H_7=\{[1,4],[2,3]\}.&&
\end{array}$$
Now observe that if the hyperplane $H_i$ already
contains the highest root $e_1-e_4$ then it cannot lead to a
maximal proper nested set. Indeed we must get a basis if we add
the highest root to a set of vectors contained in $H_i$. Thus
$H_2$, $H_3$, $H_7$ can be excluded. It remains to consider the
hyperplanes $H_1$, $H_4$, $H_5$, $H_6$.

Hyperplanes $H_1$ and $H_4$ give rise to two MPNSs each,
while $H_5$ and $H_6$ give rise to only one. So we obtain a
list of $6$ maximal nested sets (as described in Example
\ref{exam.Kn}, we identify an irreducible subset $I$ with a subset
$S$ of $[1,2,3,4]$):
$$\begin{array}{ll}
M_1=\{[1,2],[1,2,3],[1,2,3,4]\},&M_2=\{[2,3],[1,2,3],[1,2,3,4]\},\\
M_3=\{[2,3],[2,3,4],[1,2,3,4]\},&M_4=\{[3,4],[2,3,4],[1,2,3,4]\},\\
M_5=\{[1,3],[2,4],[1,2,3,4]\},  &M_6=\{[1,2],[3,4],[1,2,3,4]\}.
\end{array}$$
\end{example}

\section{Search for maximal proper nested sets adapted to a vector:
  the general case}
\label{sect.gene}

Given a vector $v$ in the cone $\CC(\CA)$, we  describe how to
search for all maximal proper nested sets belonging to
$\CP(v,\CA)$, without enumerating all MPNS.

We use as height function a linear form that is positive and
that takes different values on all elements $\alpha_i$,
and consider the total order it induces.
Let $H$ be an $\CA$-admissible hyperplane in $V$, that is,
a hyperplane spanned by a set of vectors of $\CA$.
Then the cone $\CC(\CA\cap H)$ generated by the elements of $\CA$
belonging to $H$ is a cone with non-empty interior in $H$.

We have already seen that to list all the MPNS, we have to first
list all admissible hyperplanes $H$ and then find the irreducible
components $J_1$, $J_2$, \ldots, $J_s$ of $\CA\cap H$. Then we
choose a MPNS $M_i:=\{I_i^a\}$ for $J_i$, and define $M=M_1\cup
M_2\cup M_3\cup\cdots\cup M_s\cup\{\CA\}$.

As we have seen in Example~\ref{exam.K4} we can discard some of
the hyperplanes \emph{a priori}, because they cannot lead to a
maximal proper nested set. The next lemma examines the general
situation. Let $\theta$ be the highest element in $\CA$ and $H$ a
hyperplane of $\CA$.

\begin{lemma} \label{lemm.mpns}
There exists a maximal proper nested set $M\in \CP(v,\CA)$
attached to $H$, if and only if $\theta$ does not belong to $H$
and if $v$ belongs to the cone generated by $\theta$ and $\CA\cap
H$.
\end{lemma}

\begin{proof}
The condition is necessary. Indeed $v$ must belong to the cone
generated by the elements $\theta(I_i^a)$ and $\theta$, and all
the elements $\theta(I_i^a)$ are in $\CA\cap H$. Reciprocally
consider the projection $v-\frac{\ll u,v\rr}{\ll
u,\theta\rr}\theta$, where $u$ is the equation of the hyperplane
$H$. This can be written as $v_1\oplus v_2\oplus\cdots\oplus v_s$,
where each $v_i$ is in the cone $\CC(J_i)$. Let now
$M_i\in\CP(v_i,J_i)$ be a MPNS in $J_i$. The element $v_i$ belongs
to $\CC(\theta(M_i))$. We can write
$$v=t\theta+\sum_{i=1}^s\sum_{I_i^a\in M_i} t_{i}^a\theta(I_i^a)$$
with $t_{i}^a>0$. Thus we see that the collection
$M_1\cup\cdots\cup M_s\cup\CA$ is a maximal proper nested set in
$\CP(v,\CA)$. Moreover in this way we list all elements of
$\CP(v,\CA)$.
\end{proof}

Our search for maximal proper nested sets in $\CP(v,\CA)$
will then be pursued by constructing all possible admissible
hyperplanes $H$ for which $v$ is in the convex hull of
$\CC(\CA\cap H)$ and $\theta$.
We denote by $\Hyp(v,\CA)$ the set of such $\CA$-admissible
hyperplanes.

The following easy lemma lists some obvious conditions for the set
$\Hyp(v,\CA)$. Let $u_H\in U$ be the normal vector to an
$\CA$-admissible hyperplane, meaning that $H:=\{h\in V\,|\,\ll
u_H,v\rr=0\}$.

\begin{lemma} \label{lemm.cond}
If $H\in \Hyp(v,\CA)$ then $H$ satisfies the following conditions:
\begin{enumerate}
\item $\ll u_H,\theta\rr\neq 0$.
\item $\ll u_H,v\rr\times\ll u_H,\theta\rr\geq 0$.
\end{enumerate}
\end{lemma}

Thus if a hyperplane $H$ satisfies the above conditions we
define
$$\proj_H(v)
  =v-\frac{\ll u_H,v\rr}
          {\ll u_H,\theta\rr}\theta.$$
Hence to decide if $H\in \Hyp(v,\CA)$ we simply have to test
if $\proj_H(v)$ is in the cone generated by $\CA\cap H$, which
is done by standard methods.
Our search for the hyperplanes $H\in \Hyp(v,\CA)$ will also be
considerably sped up by the following remark.

\begin{proposition} \label{prop.admi}
Let $H$ be an $\CA$-admissible hyperplane.
Let $u\in U$ be a linear form on $V$ which is non negative on
$\CA\cap H$ and on $\theta$.
If $\ll u,v\rr<0$, then $H$ is not in $\Hyp(v,\CA)$.
\end{proposition}

\begin{proof}
Indeed if $v$ was in the cone generated by $\CA\cap H$ and
$\theta$, the value of $u$ would be non negative on $v$.
\end{proof}

The point of this remark is that in classical examples of root
systems, an \emph{a priori} description of the $\CA$-admissible
hyperplanes is available, together with the defining equations of
the cone $\CC(\CA\cap H)$.
This condition will allow us to disregard right away many
$\CA$-admissible hyperplanes.

Let us summarize the scheme of the algorithm in
Figure~\ref{figu.algo.MPNS}. Recall that we have as input a vector
$v$, and as output the list of all MPNS's belonging to
$\CP(v,\CA)$.

\begin{figure}[ht]
\begin{center}
\begin{tabbing}
check if $v\in\CC(\CA)$\\
for \= each hyperplane $H$ do\\
    \> check if $v$ and $\theta$ are on the same side of $H$\\
    \> if not, then skip this hyperplane\\
    \> define the projection $\proj_H(v)$ of $v$ on $H$ along $\theta$\\
    \> check if $\proj_H(v)$ belongs to $\CC(\CA\cap H)$; if not then skip this hyperplane\\
    \> write $\CA\cap H$ as the union of its irreducible
       components $I_1\cup\cdots\cup I_k$\\
    \> write $v$ as $v_1\oplus\cdots\oplus v_k$ according
       to the previous decomposition\\
    \> for \= each $I_j$ do\\
    \>     \> compute all MPNS's for $v_j$ and $I_j$\\
    \>     \> collect all these MPNS's for $v_j$ and $I_j$\\
    \> end of loop running across $I_j$'s\\
    \> collect all MPNS's for the hyperplane $H$\\
end of loop running across $H$'s\\
return the set of all MPNS's for all hyperplanes\\
\end{tabbing}
\end{center}
\caption{Algorithm for MPNS's computation (general case)}
\label{figu.algo.MPNS}
\end{figure}

We will explain our algorithm in more details for each classical
root system (see Sections~\ref{sect.A}--\ref{sect.D}).

\section{Trees and order of poles} \label{sect.tree}

Let $M$ be a maximal nested proper set for  the system
$\CA:=\{\alpha_1,\alpha_2,\ldots,\alpha_n\}$. In our algorithms,
we will need to take an iterated residue with respect to a basis
$\overrightarrow{\theta(M)}$ of a function of the form
$\phi=\frac{P}{\prod_{i=1}^N\alpha_i}$, where $P$ is a polynomial
function on $U$. It is thus important to understand the order of
the poles of the function obtained after performing a certain
number of residues. We also  prove that the iterated residue
associated to $M$ depends only on the tree associated to $M$.

We associate to a maximal nested set $M$ a tree $T$ as follows.
Let $M=\{I_1,\ldots,I_r\}$ be a maximal nested set. The vertices
of $T$ are the elements of $M$ and the oriented edges are
determined from the reverse order relation  by  inclusion: the
ends of the tree are irreducible sets with just one element and if
$\CA$ is irreducible, the base is the set $\CA$. A subset $N$ of
$M$ will be called \emph{saturated} if it contains all elements above
elements of $N$ in the tree order. Thus if $N$ contains an element
$S$, it contains all the elements $S'$ of $M$ which are contained
in $S$.

\begin{example} \label{exam.tree.K4}
The two MNSs named $M_1$ and $M_5$ described in
Example~\ref{exam.K4} can be rewritten respectively as
$$\xymatrix{
[1,2]          &          &[1,3]&                         &[2,4]\\
[1,2,3]\ar[u]  &\mbox{and}&     &                         &\\
[1,2,3,4]\ar[u]&          &     &[1,2,3,4]\ar[uur]\ar[uul]&\\
}$$
\end{example}

Lemmas \ref{lemm.A}, \ref{lemm.B}, and \ref{lemm.D}
describe the decomposition of $\CA \cap H$ in irreducible nested
sets and lead to the following result:

\begin{proposition} \label{theo.order}
Let $T$ be the tree associated to an irreducible classical root
system. Then $T$ is a connected tree for which every vertex is
adjacent to at most two other vertices.
\end{proposition}

\begin{lemma} \label{lemm.satu}
Let
$M=\{I_1,I_2,\ldots,I_r\}$ be a MPNS. Here we have numbered our
irreducible sets such that
$\theta(I_1)<\theta(I_2)<\cdots<\theta(I_r)$.
Let  $k$  be an integer smaller than or equal to $r$.
Then the set $\{I_1,I_2,\ldots,I_k\}$ is saturated.
\end{lemma}

Indeed if two sets $I,J$ belongs to $M$ and $I\subset J$, then
$\theta(I)<\theta(J)$.

\begin{proposition}\label{reorder}
Let $M=[I_1,I_2,\ldots, I_r]$ be a maximal nested proper family.
Let $[I'_1,I'_2,\ldots,I'_r]$ be a reordering of the sequence
$[I_1,I_2,\ldots,I_r]$. We assume that this reordering is
compatible with the partial order given by inclusion: if
$I'_j\subset I'_k$ then $j<k$. Let
$$\vec{\nu}=[\theta(I_1),\theta(I_2),\ldots,\theta(I_r)]$$
and
$$\vec{\nu'}=[\theta(I'_1),\theta(I'_2),\ldots,\theta(I'_r)].$$
Then we have $\Ires_{\vec{\nu}}=\Ires_{\vec{\nu'}}.$
\end{proposition}

\begin{proof}
We prove this proposition by induction on $r$.

If $\CA$ is irreducible, then necessarily $I_r=I'_r=\CA$ and
$[I'_1,I'_2,\ldots,I'_{r-1}]$ is a reordering of the sequence
$[I_1,I_2,\ldots,I_{r-1}]$. Furthermore the families
$\{I_1,I_2,\ldots,I_{r-1}\}$ and $\{I'_1,I'_2,\ldots,I'_{r-1}\}$
are maximal proper nested sets for $\CA_0=\cup_{j=1}^{r-1}I_j$.
The set $\CA_0$ spans a codimension $1$ vector space in $V$.

To prove that $\Ires_{\vec{\nu}}=\Ires_{\vec{\nu'}}$, it suffices
to test it on basic fractions $f_\sigma$.  Let
$\sigma=\{\beta_1,\beta_2,\ldots,\beta_r\}$ be a basic subset of
$\CA$.
 By Lemma~\ref{necessary}, if $\Ires_{\vec{\nu}}f_\sigma\neq
0$, then the set $\sigma \cap \ll \CA_0\rr $ is of cardinality $r-1$,
and there exists an element of $\sigma$, say $\beta_r$,  of the
form $c\theta+\xi$ where $\xi$ belongs to $\ll \CA_0\rr $, $c$ is
a non-zero constant and $\theta$ is the highest element of $\CA$.
Let
$$\vec{\nu_0}=[\theta(I_1),\theta(I_2),\ldots,\theta(I_{r-1})]$$
and
$$\vec{\nu_0'}=[\theta(I'_1),\theta(I'_2),\ldots,\theta(I'_{r-1})].$$
Then we have
$$\Ires_{\vec{\nu}}f_\sigma
  =\frac{1}{c}\,\Ires_{\vec{\nu_0}}f_{\sigma\cap\ll \CA_0\rr }$$
and
$$\Ires_{\vec{\nu'}}f_\sigma
  =\frac{1}{c}\,\Ires_{\vec{\nu'_0}}f_{\sigma\cap\ll \CA_0\rr }.$$
We conclude by induction.

When $\CA$ is not irreducible, we write $\CA=\cup_{a=1}^s J_a$
where $J_a$ are irreducibles.
We have $V=\oplus_{a=1}^s \ll J_a\rr $.
Every basic subset $\sigma$ of $\CA$ is the union of basic
subsets for the irreducible sets $J_a$.
Define
$$\vec{\nu}_a=[\theta(I_a^{i_1}),\theta(I_a^{i_2}),\ldots,\theta(J_a)]$$
where $[I_a^{i_1},I_a^{i_2},\ldots,J_a]$ is the subsequence of
irreducible sets contained in $J_a$ extracted (with conserving
order) from the sequence $[I_1,I_2,\ldots,I_r]$.
Similarly let
$$\vec{\nu'}_a=[\theta({I'_a}^{i_1}),\theta({I'_a}^{i_2}),\ldots,\theta(J_a)]$$
where $[{I'_a}^{i_1},{I'_a}^{i_2},\ldots,I'_a]$ is the subsequence of
irreducible sets contained in $J_a$ extracted from the sequence
$[I'_1,I'_2,\ldots,I'_r]$. Then, as the calculation takes place
with respect to independent variables, we have
\begin{eqnarray*}
\Ires_{\vec{\nu}}(f_\sigma)
  &=&\prod_{a=1}^s(\Ires_{\vec{\nu}_a}f_{\sigma\cap\ll J_a\rr }),\\
\Ires_{\vec{\nu'}}(f_\sigma)
  &=&\prod_{a=1}^s(\Ires_{\vec{\nu'}_a}f_{\sigma\cap\ll J_a\rr }).
\end{eqnarray*}
Each of the  vector space $\ll J_a\rr $ is of dimension less than
$r$, so that by induction hypothesis
$\Ires_{\vec{\nu}_a}=\Ires_{\vec{\nu'}_a}$. This concludes the
proof.
\end{proof}


Let us  now consider partial iterated residues. To a set $\nu$ of
elements of $\CA$, we associate the vector space
$$H_\nu
  :=\{u\in U\,|\,\ll\alpha,u\rr=0\mbox{ for all }\alpha\in\nu\}.$$
A linear function $\alpha\in \CA$ produces a linear function on
$H_\nu$ by restriction.
If $\vec{\nu}:=[\alpha_1,\alpha_2,\ldots,\alpha_k]$ is a sequence
of elements of $\CA$, the partial iterated residue
$$\Ires_{\vec{\nu}}\phi:=\res_{\alpha_k=0}\cdots\res_{\alpha_1=0}\phi$$
associates to a rational function $\phi$ in $R_{\CA}$ a rational
function on $H_\nu$ of the form

$$\frac{G}{\prod_{i=1,\ldots,n\,;\,\overline{\alpha_i}\neq 0}
          {\overline \alpha_i}^{n_i}},$$
where $G$ is
a polynomial function on $H_\nu$ and $\overline\alpha$ is the
restriction of $\alpha$ to $H_\nu$.
 Let $M$ be a MPNS and consider the tree associated to $M$.
Given a saturated subset $S$ of $M$, we can  define the iterated
residue with respect to this saturated set: we choose any order
$S:=[I_1,I_2,\ldots,I_k]$ on $S$ compatible with the inclusion
relation and define $\Ires_S:=\Ires_{\vec{\nu}}$ with
$\vec{\nu}=[\theta(I_1),\theta(I_2),\ldots,\theta(I_k)]$. With the
same proof as for Proposition~\ref{reorder}, this partial residue
depends only on the set $S$.
 We denote by
$H_S$ the intersection of the kernels of the elements $\alpha$ for
$\alpha \in S$. It is also the intersection of the kernels of the
elements $\theta(I_k)$, as the set $\nu$ is a basic sequence in
$S$.


 Let $\phi$ be a function in $R_{\CA}$ of the form
$$\phi=\frac{P}{\prod_{i=1}^n\alpha_i}.$$

Let $M$ be a MPNS and $J_1,J_2,\ldots, J_s$ be elements of $M$. We
consider the saturated subset $S$ of $M$ consisting of the
elements of the tree \emph{strictly above} $J_1,J_2,\ldots,J_s$.
The iterated residue $\Ires_S\phi$ is a function on $H_S$. Denote
by $u_a$ the restriction of the function $\theta(J_a)$ to $H_S$.

\begin{proposition}\label{prop.order}
The pole of the linear function $u_a$ in the iterated residue
$\Ires_S\phi$ is of order less than or equal to $|J_a|-\dim\ll
J_a\rr+1$.
\end{proposition}

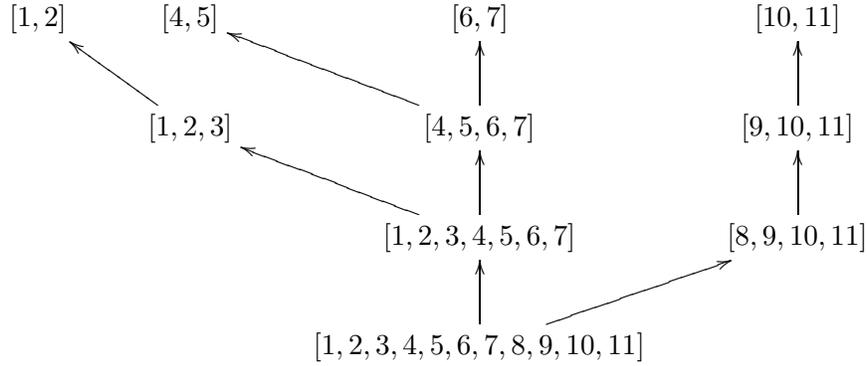
\begin{figure}[ht]
$\xymatrix{
[1,2]\ar@{<-}[dr]&  [4,5]\ar@{<-}[dr]&          [6,7]\ar@{<-}[d]&     [10,11]\ar@{<-}[d]\\
                 &[1,2,3]\ar@{<-}[dr]&      [4,5,6,7]\ar@{<-}[d]&   [9,10,11]\ar@{<-}[d]\\
                 &                   &[1,2,3,4,5,6,7]\ar@{<-}[d]&[8,9,10,11]\ar@{<-}[dl]\\
                 &                   &[1,2,3,4,5,6,7,8,9,10,11]&
}$
\caption{Irreducible components of a MPNS in $A_{10}$}
\label{figu.irred}
\end{figure}

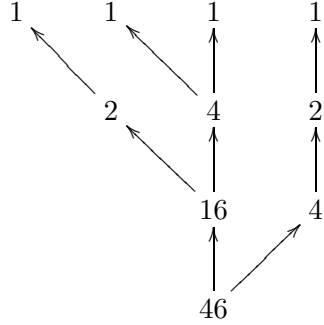
\begin{figure}[ht]
$\xymatrix{
1\ar@{<-}[dr] & 1\ar@{<-}[dr] &  1\ar@{<-}[d] & 1\ar@{<-}[d]\\
              & 2\ar@{<-}[dr] &  4\ar@{<-}[d] & 2\ar@{<-}[d]\\
              &               & 16\ar@{<-}[d] & 4\ar@{<-}[dl]\\
              &               & 46            &
}$
\caption{Order of nodes in the tree represented in
  Figure~\ref{figu.irred}, according to
  Proposition~\ref{prop.order}}
\label{figu.tree}
\end{figure}

See Figures~\ref{figu.irred} and~\ref{figu.tree} for an
application of the proposition.

\begin{proof}
Choose a vector space $E$ such that
$$V=\ll J_1\rr\oplus\cdots\oplus\ll J_s\rr\oplus E.$$
Let $B=\bigcup_{a=1}^sJ_a$ and $C=\CA\setminus B$.
Write $C:=\{\beta_1,\beta_2,\ldots,\beta_q\}$
and
$$\phi=P\times\phi_1\times\phi_2\times\cdots\times\phi_s\times Q$$
with $\phi_a=\frac{1}{\prod_{\alpha\in J_a}\alpha}$ and
$Q=\frac{1}{\prod_{j=1}^q\beta_j}$.

For $\beta_j\in C$, we write
$\beta_j=\sum_{i=1}^s\beta_j^i+\gamma_j$ with $\beta_j^i\in\ll
J_i\rr$ and $\gamma_j\in E$. The element $\gamma_j$ is necessarily
non zero, as the set $B$ is complete. Thus we write
$$\frac{1}{\beta_j}
  =\frac{1}
        {\gamma_j\left(1+\frac{\sum_{i=1}^s\beta_j^i}{\gamma_j}\right)}$$
and the iterated residue is by definition
$$\Ires_S(\phi)
  =\Ires_S\left(P\times(\phi_1\cdots\phi_s)
                 \times\prod_{j=1}^q\frac{1}{\gamma_j}
                 \sum_{k=0}^{\infty}
                    \left((-1)^k\frac{\sum_{i=1}^s\beta_j^i}{\gamma_j}
                    \right)^k
          \right).$$
Here, when taking the residue, the elements $\gamma_j$ are considered
as constants and this sum is finite.

Consider the subset $M_a$ of elements of $M$ contained in $J_a$.
This is a MPNS for the set $J_a$. Let $J_a^+$ be the saturated
subset of $M_a$ consisting of all elements of $M_a$ different from
$J_a$. Then $J_a^+$ has $\dim\ll J_a\rr -1$ elements. If
$g=\frac{P_a}{\prod_{\alpha\in J_a}\alpha^{n_\alpha}}$, the
iterated residue $\Ires_{J_a^+}g$ is a Laurent  polynomial in
$u_a$.

Now $\Ires_S\phi$  is a sum of products of residues of the form
$\Ires_{J_a^+}g_a$ where $g_a=\frac{P_a}{\prod_{\alpha\in
J_a}\alpha}$ and $P_a$ is a polynomial. Thus we obtain a Laurent
polynomial in $u_a, a=1,\ldots, s$ (with coefficients  rational
functions on the vector space $E^*$).
 Now the homogeneous degree
of $g_a$ is greater than or equal to $-|J_a|$. The number of residues
we are taking is equal to $\dim \ll J_a\rr - 1$.
So we obtain a function of $u_a$ of homogeneous degree greater than or
equal to $-|J_a|+\dim\ll J_a\rr-1$. This means that the pole in $u_a$
is of order less than or equal to $|J_a|-\dim\ll J_a\rr+1$.
\end{proof}

Let us consider the MPNS whose tree representation is given
by Figure~\ref{figu.irred}.
The orders of the poles of its nodes are given in
Figure~\ref{figu.tree}.

\begin{remark} \label{remark.ordertree}
In our program for calculating iterated residues for root systems
of type $A_r$, we will reorder roots according to the tree order: we
take the residue first with respect to the elements $\theta(I_k)$
appearing at the end of the tree in arbitrary order, and we remove
these variables. Then we take the variables appearing at the end
of the tree when we have removed these irreducible sets. Here an
irreducible set $I$ is indexed by a subset $S$ of
$\{1,2,\ldots,r+1\}$. A subset $S$ of cardinality $2$, for example
$[1,3]$, corresponds to the irreducible set with one element (here
$e_1-e_3$). Thus given a MNS $M$ represented as
$M=\{S_1,S_2,\ldots,S_r\}$ we will first take the residues with
respect to the roots $\theta(I_k)$, for sets $S_k$ of cardinality
$2$, in arbitrary order, then with respect to irreducible sets
associated to sets $S_k$ of cardinality $3$, \etc The procedure of
ordering roots  coming from a MNS $M=\{S_j\}$ according to the
cardinality of the set $S_k$ is called ${\tt OrderThetas}$.
Furthermore we will at the same time keep track of the order of
the pole for calculating an iterated residue of a function
$\phi=P/\prod_{i=1}^N\alpha_i$ in the procedure ${\tt
FormalPathAwithOrders}$.
\end{remark}


\section{Volume and partition function for the system $A_{n-1}$}
\label{sect.A}

\subsection{The formulae to be implemented}
\label{sect.A.form}

Let $E$ be an $n$-dimensional vector space with basis $e_i$
($i=1$, \ldots, $n$) and consider the set
$$\CK_n=\{e_i-e_j\,|\,1\leq i<j\leq n\}.$$
These are the positive roots for a system of type $A_{n-1}$.
The number of elements in $\CK_n$ is $N=n(n-1)/2$.
Note that $\CK_n$ is also the set of vectors in a
complete graph with $n$ nodes.

We let $V$ be the vector space generated by the elements in $\CK_n$.
Then $V$ has dimension $n-1$ and it is defined by:
$$V=\left\{v=\sum_{i=1}^n v_ie_i\in E\,\Big|\,
           \sum_{i=1}^n v_i=0\right\}.$$
In our procedures, a vector $v$ of length $n$ such that
$\sum_{i=1}^n v_i=0$ will be called an \emph{$A$-vector} and
written as $v=[v_1,v_2,\ldots,v_n]$.
The lattice spanned by $\CK_n$ is simply
$$V_\Z=\left\{h=\sum_{i=1}^nh_ie_i\in\Z^n\,\Big|\,\sum_{i=1}^nh_i=0\right\}.$$
It is well known and easy to prove that $\CK_n$ is unimodular. The
cone $\CC(\CK_n)$ generated by $\CK_n$ is simplicial with
generators the $n-1$ simple roots
$e_1-e_2,e_2-e_3,\ldots,e_{n-1}-e_n$. This cone is described as:
$$\CC(\CK_n)=\{\mbox{$A$-vector }v=[v_1,v_2,\ldots,v_n]\,|\,
    v_1+v_2+\cdots+v_i\geq 0\mbox{ for all $i$}\}.$$
Keep in mind that our vector $v$ satisfies the condition
$$v_1+v_2+\cdots+v_{n-1}+v_n=0.$$

We choose on $V$ the measure $dh$ determined by $V_\Z$.
Let $v$ be in
the cone $\CC(\CK_n)$. We are interested to compute the volume
$\vol_{\Z,\CK_n}(v)$ of the polytope
$$\Pi_{\CK_n}(v)
  =\left\{(x_\alpha)_\alpha\in\R^N\,\Big|\,x\geq 0,
          \sum_{\alpha\in\CK_n}x_\alpha\alpha=v\right\}.$$
If $h$ is a point in $V$ with integral coordinates then
we are also interested in computing the number $N_{\CK_n}(h)$
of integral points in $\Pi_{\CK_n}(h)$.

We apply the formulae of Theorem~\ref{theo.main}. Since $\CK_n$ is
unimodular, the set $F$ can be taken as $F:=\{0\}$
(Remark~\ref{rema.unim}).

Since $V$ is contained in $E$, then we have a canonical map
$E^*\lra V^*$ given by restriction.
Define $U=V^*$ as in the general setting.
We identify $U$ with $\R^{n-1}$ by sending $u\in\R^{n-1}$ to
$u=\sum_{i=1}^{n-1}u_i e^i\in E^*$, where $e^i$ is the dual
basis to $e_i$.
Thus the root $e_i-e_j$ ($1\leq i<j<n$) produces the linear
function $u_i-u_j$ on $U$, while the root $e_i-e_n$ produces
the linear function $u_i$.


\begin{definition} \label{defi.func}
Let $v=\sum_{i=1}^n v_ie_i\in V$ be a vector with real
coordinates.
Let $h=\sum_{i=1}^nh_ie_i\in V$ be a vector with integral
coordinates.
Then for $u\in U$ define:
\begin{itemize}
\item
  $\displaystyle J_{A}(v)(u)=
    \frac{e^{\sum_{i=1}^{n-1}u_iv_i}}
   {\prod_{i=1}^{n-1}u_i\prod_{1\leq i<j\leq n-1}(u_i-u_j)}$
\item
  $\displaystyle \CF_{A}(h)(u)=
    \frac{\prod_{i=1}^{n-1}(1+u_i)^{h_i+n-1-i}}
         {\prod_{i=1}^{n-1}u_i\prod_{1\leq i<j\leq n-1}(u_i-u_j)}$ .
\end{itemize}
\end{definition}

\begin{theorem} \label{theo.A}
Let $\gc$ be a chamber of $\CC(\CK_n)$.
\begin{itemize}
\item For $v\in\overline\gc$, we have
  $$\vol_{\Z,\CK_n}(v)=\JK_\gc(J_{A}(v)).$$
\item For $h\in\Z^n\cap\overline\gc$, we have
  $$N_{\CK_n}(h)=\JK_\gc \left(\CF_{A}(h) \right) .$$
\end{itemize}
\end{theorem}

\begin{proof}
The first assertion is the general formula.

The function $F(0,h)(u)=e^{\ll h,u\rr}/
  \prod_{\alpha\in\CA}(1-e^{-\ll\alpha,u\rr})$
for the system $\CK_n$ is
$$F(0,h)(u)
  =\frac{e^{\sum_{i=1}^{n-1}u_iv_i}}
        {\prod_{i=1}^{n-1}(1-e^{-u_i})
         \prod_{1\leq i<j\leq n-1}(1-e^{-(u_i-u_j)})}.$$
Note that the change of variable $1+z_i=e^{u_i}$ preserves the
hyperplanes $u_i=0$ and $u_i=u_j$.
After the change of variable, we get
\begin{equation} \label{equa.para}
F(0,h)(u)=\frac{\prod_{i=1}^{n-1}(1+z_i)^{h_i+n-i}}
     {\prod_{1\leq i<j\leq n-1}(z_i-z_j)
      \times
      \prod_{i=1}^{n-1}z_i}.
\end{equation}
But $z_i=e^{u_i}-1$ leads to $dz_i=e^{u_i}du_i=(1+z_i)du_i$
and hence we obtain the desired exponent $h_i+n-i-1$ thanks
to the formula involving Jacobians in
Proposition~\ref{prop.chng}.
\end{proof}

In order to implement these formulae, we first have to describe
the set $\CP(v,\CK_n)$ (Section~\ref{sect.A.mpns}), then calculate
the iterated residue formulae associated to these paths
(Section~\ref{sect.A.ires}). Below we explain how these
computations fit together to get a global procedure for the Kostant
partition function for $A_{n-1}$ (Section~\ref{sect.A.glob}). As a
short digression, we will explain how we adapted our program to
deal with formal parameters (Section~\ref{sect.A.para}).
\subsection{The search for maximal proper nested sets adapted
  to a vector}
\label{sect.A.mpns}

We now look for maximal proper nested sets adapted to a vector
following the general method as outlined in
Figure~\ref{figu.algo.MPNS}: we will begin by listing all possible
$\CK_n$-admissible hyperplanes.
The usual height function is
$$\h(v)=\sum_{i=1}^{n-1}(n-i)v_i$$
which takes the value $1$ on all the simple roots, and hence the
value $j-i$ on $e_i-e_j$.
We deform $\h$ slightly in order to have a function taking
different values on all roots:
If two elements $e_i-e_j$ and $e_k-e_\ell$ are such that
$j-i=\ell-k$, we decide that $\h(e_i-e_j)<\h(e_k-e_\ell)$ if $i<k$.

If $P$ is a proper subset of $\{1,2,\ldots,n\}$ and $v$ is an
$A$-vector, we denote by $\ll u_P,v\rr$ the linear form
$\sum_{i\in P}v_i$, and by $H_P$ the hyperplane
$$H_P:=\{v\in V\,|\,\ll u_P,v\rr=0\}.$$

We will see shortly that all $\CK_n$-admissible hyperplanes are
obtained in this way, that is giving a proper subset $P$ of
$\{1,\ldots,n\}$. Observe that the hyperplane $H_P$ is equal to
the hyperplane $H_Q$ determined by the complement $Q$ of $P$.
We denote
$$\CK(P):=
  \{e_i-e_j\,|\,\,1\leq i<j\leq n\,;\,i,j\in P\}\subset\CK_n.$$
Note that $\CK(P)$ is the positive system
$A_{|P|-1}$, where the positivity is induced by the
lexicographic order.

\begin{lemma} \label{lemm.A}
\begin{itemize}
\item The hyperplane $H_P$ is a $\CK_n$-admissible hyperplane.
\item The set $\CK_n\cap H_P$ is the union of $\CK(P)$ and
  $\CK(Q)$, where $Q$ is the complement of $P$ in
  $\{1,2,\ldots,n\}$.
\item Every $\CK_n$-admissible hyperplane is of this form.
\end{itemize}
\end{lemma}

\begin{proof}
The first two assumptions are easy to see.
We prove the third by induction on $n$, the case $n=2$ being
trivial.
Let $H$ be a $\CK_n$-admissible hyperplane.
Let $\alpha$ be a root in $H$.
Renumbering the roots, we may assume that $\alpha=e_{n-1}-e_n$.
The map $q$ sending $e_i$ to $e_i$ if $i<n$ and $e_n$ to
$e_{n-1}$ sends the set $\CK_n\setminus\{\alpha\}$
to $\CK_{n-1}$.
The space $H/\R\alpha$ becomes a $\CK_{n-1}$-admissible
hyperplane.
It is thus determined by a subset $P'$ of $\{1,2,\ldots,n-1\}$.
If $P'$ does not contain $n-1$, the hyperplane $H$ is equal to
the hyperplane determined by the subset $P'$ of
$\{1,2,\ldots,n-1\}$.
If $P'$ does contain $n-1$, then the hyperplane $H$ is equal to
the hyperplane determined by $P=P'\cup\{n\}$.
\end{proof}

We now proceed to the detailed description of our algorithm.
Recall our description of an $A$-vector as an array
$v=[v_1,v_2,\ldots,v_n]$ with $\sum_{i=1}^n v_i=0$.
Referring to Figure~\ref{figu.algo.MPNS}
we need to check if the vector is in the cone $\CC(\CK_n)$, that
is, $\sum_{j=1}^iv_j\geq 0$ for $1\leq i\leq n-1$.
This is done by using the procedure {\tt CheckVector(v)},
which gives an answer true or false.

For the system $\CK_n$ the highest root $\theta$ is equal to
$$\theta=[1,0,\ldots,0,-1]\in\R^n$$
and computed with the procedure ${\tt theta}$.

At this point, we need to list all hyperplanes that are in
$\Hyp(v,\CK_n)$. This is done in the procedure ${\tt TwoSets(v)}$,
that we are about to describe. As explained in Lemma~\ref{lemm.A},
each hyperplane is determined by an equation $\sum_{i\in I}
a_i=0$. It therefore produces a set of two lists $P$, $Q$, where
$P=[i\in I]$ and $Q=[i\notin I\,|\,1\leq i\leq n]$. Note that $P$
and $Q$ are sorted. To verify that such a hyperplane is in
$\Hyp(v,\CK_n)$, we need to test if $\ll u_P,\theta\rr$ is not zero
($\theta$ is not in the hyperplane) and if $\ll u_P,v\rr\times\ll
u_P,\theta\rr$ is non-negative ($\theta$ and $v$ are on the same
side of the hyperplane).

Furthermore, the procedure ${\tt ProjH(v,H)}$ constructs the
vector
$$\proj_H(v)
  =v-\frac{\ll u_P,v\rr}
          {\ll u_P,\theta\rr}\theta$$
that we represent as $\{[v_1,P],[v_2,Q]\}$.
Each of the vectors $v_1$, $v_2$ is an $A$-vector (sum of
coordinates equal to zero).
So the last condition for $H$ being in $\Hyp(v,\CK_n)$ is that
$v_1\in\CC(\CK(P))$ and $v_2\in\CC(\CK(Q))$.

Hence a hyperplane $H$ is in $\Hyp(v,\CK_n)$ if it satisfies the
series of conditions:
\begin{eqnarray*}
\ll u_P,\theta\rr\neq 0
  &\mbox{with}&{\tt Hvalue(theta(n),P)\neq 0},\\
\ll u_P,v\rr\times\ll u_P,\theta\rr\geq 0
  &\mbox{with}&{\tt CheckSide(v,P)=true},\\
v_1\in\CC(\CK_{|P|-1})
  &\mbox{with}&{\tt CheckVector(v_1)=true},\\
v_2\in\CC(\CK_{|Q|-1})
  &\mbox{with}&{\tt CheckVector(v_2)=true}.\\
\end{eqnarray*}
The procedure ${\tt CheckList(v,H)}$ implements all
these sub-routines. It is used in the procedure ${\tt
TwoSets(v)}$, computing all elements of $\Hyp(v,\CK_n)$.
We combine ${\tt TwoSets}$ with a procedure named ${\tt
TwoVector}$ to finally get the procedure ${\tt TwoVectors(v)}$
determining all hyperplanes in  $\Hyp(v,\CK_n)$ and projections of
$v$ on these hyperplanes.

We now have to perform the next step of our algorithm.
Let $\{[v_1,K_1],[v_2,K_2]\}$ be the output of
${\tt TwoVectors(v)}$.
Then we construct the MNSs for $[v_1,K_1]$ and $[v_2,K_2]$,
and go on recursively until the procedure stops.
These iterated steps are done by the procedure ${\tt Splits}$.

Finally the procedure ${\tt MNSs(v)}$, computing all MNSs for a
given vector $v$, works as follows.
We begin by building the first seed of MNSs with the procedure
${\tt MNS1}$, containing the regularization of the result of
${\tt TwoVectors}$.
We then call repeatedly the procedure ${\tt AllNewMNSs}$,
which performs the regularization of the output of
${\tt Splits}$.

\subsection{Residues associated to maximal proper nested sets}
\label{sect.A.ires}

An element $M$ in $\CP(v,\CK_n)$ is represented as a collection
$M=\{K_1,K_2,\ldots,K_{n-1}\}$ of $(n-1)$  subsets of
$[1,2,\ldots,n]$.
As we have said in Remark~\ref{remark.ordertree}, given a maximal
proper nested set $M:=\{K_1,K_2,\ldots,K_{n-1}\}$ we associate
to it an ordered basis $\overrightarrow{\theta(M)}$ of $V$
(procedure ${\tt OrderThetas}$).
 If $p=[\alpha_1,\alpha_2,\ldots,\alpha_{n-1}]$ is
the list of roots singled out by our procedure, then  $\alpha_1$
is an element associated to a set $K_i$ of cardinality $2$ and
$\alpha_{n-1}=\theta$. We identify the root $e_i-e_n$ to the
linear function $z_i$  on $\C^{n-1}$ and the root $e_i-e_j$ to
$z_i-z_j$.

Let $h$ be a $A$-vector with integral coordinates. Let us consider
the Kostant function
$$\CF_A(h)(z_1,z_2,\ldots,z_{n-1})=\frac{\prod_{i=1}^{n-1}(1+z_i)^{h_i+n-1-i}}
     {\prod_{1\leq i<j\leq n-1}(z_i-z_j)
      \times
      \prod_{i=1}^{n-1}z_i}$$
(Definition~\ref{defi.func}). To compute $N_{\CK_n}(h)$,
we will have to compute
$$\res_M\phi:=
  \res_{\alpha_{n-1}=0}\res_{\alpha_2=0}\cdots\res_{\alpha_1=0}\phi$$
with $\phi=\CF_A(h)$.
Using Proposition~\ref{prop.order}, we know in advance the order
of the pole in $\alpha_k=0$ of the function obtained after taking
the first $(k-1)$ residues.
These orders are recorded in the procedure
${\tt FormalPathAwithOrders}$.


If $\alpha_1=z_i-z_j$, we can replace --- after taking
the residue at $z_i=z_j$ --- the variable $z_i$ by the variable
$z_j$ in all the other roots. Thus we get rid of the variable
$z_i$. The procedure ${\tt NewR}$ produces the ordered path
resulting from all these substitutions.

Recursively, we will have  to compute the residue at $\rt=0$ of an
expression

\begin{equation} \label{equa.resi}
f=\frac{A(z_i,i\in L)}
         {\prod_{i,j\in L; i<j}(z_i-z_j)^{m_{i,j}}
          \prod_{i\in L}z_i^{m_i}
         },
\end{equation}
where $L$ is a list of indices taken in
$\{1,\ldots,n-1\}$. Denote by $\maxi$ the order $m_{i_0,j_0}$ of
the root $\rt$ (the exponent $\maxi$ is recorded in the procedure
${\tt FormalPathAwithOrders}$). Note that computing the residue
is exactly the same as computing the coefficient of $z$ of degree
$\maxi-1$ of the expansion of $f\times(\rt)^{\maxi}$ at
$z_{i_0}=z+z_{j_0}$. Let us describe in detail the procedure ${\tt
ComputeRes}$, performing this task.

For $j\in L\setminus\{i_0,j_0\}$, let
et $e_j=m_{i_0,j}$ if $i_0<j$ and $e_j=m_{j,i_0}$ if $i_0>j$.
Then $f\times(\rt)^{\maxi}$ can be written as
$g\times h\times R$ where

\begin{eqnarray} \label{equa.g}
g&=&\frac{A(z_i,i\in L)}
         {\prod_{i,j\in L; i,j\neq i_0;i<j}(z_i-z_j)^{m_{i,j}}}
    \times
    \frac{1}{\prod_{i\in L;i\neq i_0}z_i^{m_i}},\\ \label{equa.h}
h&=&\frac{1}{\prod_{j\in L;j\neq i_0,j_0}(z_j-z_{i_0})^{e_j}}
    \times
    \frac{1}{z_{i_0}^{e_{0}}},\\ \label{equa.R}
R&=&(-1)^{\sum_{j>i_0}e_j}.
\end{eqnarray}

Let $z_{i_0}=z+z_{j_0}$. So to get the desired residue, we need to
calculate the expansion of $g$ and $h$ at $z=0$. More precisely if
$g=\sum_{i=0}^{\maxi-1}g_iz^i$ and $h=\sum_{i=0}^{\maxi-1}h_iz^i$,
then the coefficient of degree $\maxi-1$ of $f\times(\rt)^{\maxi}$
is simply $R\times\sum_{i=0}^{\maxi-1}g_ih_{\maxi-1-i}$. Let us
describe how the procedure ${\tt ComputeRes}$ performs this task.

Rewrite the fraction $h$ defined in
Equation~\eqref{equa.h} as $B\times\tilde{h}$, where

\begin{eqnarray*} 
B&=&\frac{1}{\prod_{j\in L;j\neq i_0,j_0}(z_j-z_{j_0})^{e_j}}
    \times
    \frac{1}{z_{j_0}^{e_0}},\\ 
\tilde{h}&=&\frac{1}
                 {\prod_{j\in L;j\neq i_0,j_0}
                  \left(1-\frac{z}{z_j-z_{j_0}}\right)^{e_j}}
            \times
            \frac{1}{\left(1+\frac{z}{z_{j_0}}
                     \right)^{e_0}}.
\end{eqnarray*}
Consequently, to expand $h$ as a function of $z$, we only
need to expand $\tilde{h}$. This is done in the procedure ${\tt
CoeffBin}$ using the binomial coefficients.
In the procedure ${\tt CoeffFun}$ we calculate the expansion
at $\rt=0$ of the fraction

\begin{eqnarray} \label{equa.cfun}
& & f\times(\rt)^{\maxi}
     \times z_{i_0}^{e_0}
     \times\prod_{j\in L; j\neq i_0,j_0}(z_{i_0}-z_j)^{e_j}\\ \nonumber
&=& g\times R\times(-1)^{\sum_{j\neq i_0,j_0}e_j}.
\end{eqnarray}
Finally the procedure ${\tt ComputeRes}$ performs the
sum over $i$ ranging from $0$ to $\maxi$ of

\begin{center}
\begin{tabular}{cl}
        &($S:=(-1)^{\sum_{j\neq i_0,j_0}e_j}$)
         $\times$
         $B$\\
$\times$&(the component of degree $i$ of ${\tt CoeffFun}$)\\
$\times$&(the component of degree $\maxi-1-i$ of ${\tt CoeffBin}$).
\end{tabular}
\end{center}
Rewrite this as the sum over $i$ of

\begin{center}
\begin{tabular}{cl}
        &$R$\\
$\times$&(the component of degree $i$ of ${\tt CoeffFun}$)
         $\times$ $R$ $\times$ $S$\\
$\times$&$B$ $\times$
         (the component of degree $\maxi-1-i$ of ${\tt CoeffBin}$),
\end{tabular}
\end{center}
or, equivalently, as
$\sum_{i=0}^{\maxi-1}R\times g_i\times h_{\maxi-1-i}$:
this is exactly the desired coefficient.

\begin{remark} \label{rema.sr}
For residues along roots of type $z_{i_0}$ instead of $\rt$ the
procedure ${\tt ComputeRes}$ also calls procedures ${\tt
srCoeffFun}$ and ${\tt srCoeffBin}$, similar to ${\tt CoeffFun}$
and ${\tt CoeffBin}$.
\end{remark}

\subsection{The procedure ${\tt MNS\_KostantA}$}
\label{sect.A.glob}

We finish the section dedicated to $A_{n-1}$ by giving the global
outline of the procedure ${\tt MNS\_KostantA(v)}$ computing the
Kostant partition number of a vector $v$ lying in the root
lattice. We begin by slightly deforming $v$ so that it lies on no
admissible hyperplanes, with the command ${\tt
v':=DefVector(v,n)}$. We compute all MPNSs for $v'$ with the
procedure ${\tt MNSs(v')}$.

Given such a MPNS $M=\{S_k\}$, we extract the
highest roots of its irreducible components with the call ${\tt
R:=ThetaMNS(M)}$. We obtain a set $R$  where each element of $R$
is a root represented as $[i,j]$ together with the cardinality of the
set $S_k$ it comes from.  We then transform this set $R$ into a
path $p$  keeping track of order of poles  by setting ${\tt
p:=FormalPathAwithOrders(R)}$.

Finally we compute the residue associated to this path with the
command ${\tt OneIteratedResidue(p,v,n)}$. Summing all these
residues over the set of MNSs, we obtain thanks to
Theorem~\ref{theo.A}, the desired partition number for $v$.

Let us describe in detail the procedure ${\tt
OneIteratedResidue(p,v,n)}$ computing the iterated residue along a
path $p$ for a vector $v$ lying in the root lattice for $A_{n-1}$.
We first compute the Kostant fraction (second item of
Definition~\ref{defi.func}, procedure ${\tt KostantFunctionA}$).
Then we replace in the path all roots $z_i-z_n$ by $z_i$ (with
${\tt Kpath}$). We also build  a upper bound for the orders of the
roots $m_{ij}$ (with ${\tt Multiplicity})$. Keep in mind that the
exponent $\maxi$ needed in the residue calculation is computed
\emph{a priori}, and our computation seems quite optimal.
Then we compute iteratively the residues, using the
procedure ${\tt ComputeRes}$ (Section~\ref{sect.A.ires}).
Note that at each step we have to update the list of orders
(with ${\tt MultRoots}$) and the list of remaining variables
(with ${\tt ListOfVariables}$).

\subsection{Parametrized version of the algorithm}
\label{sect.A.para}

Our algorithm can work with formal parameters, only needing slight
modifications of procedures.  We are then able to compute directly
the polynomial $h\mapsto N_{\CK_n}(h)$ giving the number of
integral points in the polytope $\Pi_{\CK_n}(h)$, on the chamber
determined by $h$ (this chamber is easily computed).
As a consequence we can easily get the Ehrhart polynomial
$t\mapsto N_{\CK_n}(th_1,\ldots,th_n)$ of the
polytope. See~\cite{Ehr74}.

Now let us outline how this modified program works. Given an
element $h=(h_1,\ldots,h_n)$ of the root lattice for $A_{n-1}$, we
want to compute the Kostant partition function for the vector
$(h_1,\ldots,h_n)$, when $h$ varies in $\Z^n\cap \overline \gc$ .

Recall that we have to perform the residue at $\rt=0$ of the
fraction defined in Equations~\eqref{equa.para}
and~\eqref{equa.resi}.
Note that the numerator of the fraction $\CF_A(h)$ contains terms
of the form $(1+z_i)^{h_i+n-1-i}$ that we must formally expand.
For any parameter $b$, we write at $z_{i_0}=z+z_{j_0}$,
$$(1+z_{i_0})^b=(1+z_{j_0})^b\left(1+\sum_{j=1}^{\maxi-1}{\binome bj}
                                     \left(\frac{z}{1+z_{j_0}}
                                     \right)^j\right)+O(z^{\maxi}).$$
Hence, for any parameter $b$, the procedure ${\tt CoeffFun}$
(see Equation~\eqref{equa.cfun}) now computes the expansion at
$\rt=0$ of the fraction

\begin{eqnarray*}
f&\times&\left(z_{i_0}^{e_0}\times
                \prod_{j\in L; j\neq i_0,j_0}(z_{i_0}-z_j)^{e_j}
          \right)\\
&\times&\left(1+\sum_{j=1}^{\maxi-1}{\binome bj}
                                     \left(\frac{\rt}{1+z_{j_0}}
                                     \right)^j
        \right)
        \times(1+z_{j_0})^{\maxi-1}.
\end{eqnarray*}
For the residue along a root of type $z_{i_0}$ instead of $\rt$,
the procedure ${\tt srCoeffFun}$ has been modified in a similar
way.


\section{The type $B_n$}
\label{sect.B}

\subsection{The formulae to be implemented}
\label{sect.B.form}

Consider a vector space $V$ with basis $e_1$, $e_2$, \ldots,
$e_n$. We choose on $V$ the standard Lebesgue measure $dh$.
Let

$$\CB_n=\{e_i\,|\,1\leq i\leq n\}
  \cup\{e_i-e_j\,|\,1\leq i<j\leq n\}
  \cup\{e_i+e_j\,|\,1\leq i<j\leq n\}.$$
Then $\CB_n$ is a positive roots system of type $B_n$ and
generates $V$.
The number of elements in $\CB_n$ is $N=n^2$.
We denote by $U$ the dual of $V$.
The lattice $V_\Z$ generated by roots is equal to $\Z^n$,
so the constant $\vol \left( V/V_\Z, dh \right) = 1$.

The cone $\CC(\CB_n)$ is simplicial and spanned by the $n$ simple
roots $e_1-e_2$, $e_2-e_3$, \ldots, $e_{n-1}-e_n$, $e_n$.
A vector $v=[v_1, v_2,\ldots,v_n]$ is in $\CC(\CB_n)$ if and only
if it satisfies the inequations $v_1+\cdots+v_i\geq 0$ for all
$i=1$, \ldots, $n$.

Let $v$ be in the cone $\CC(\CB_n)$.
Consider the polytope

$$\Pi_{\CB_n}(v)
  =\left\{(x_\alpha)_\alpha\geq 0\,\Big|\,
    \sum_{\alpha\in\CB_n}x_\alpha\alpha=v\right\}.$$
If $h$ is a point in $V$ with integral coordinates, we are
interested in computing the number $N_{\CB_n}(h)$ of
integral points in $\Pi_{\CB_n}(h)$.

Let $U_\Z$ be the lattice dual to $V_\Z$.
We identify the torus $T=U/U_\Z=\R^n/\Z^n$ to $(S^1)^n$ by

$$(u_1,u_2,\ldots, u_n)\mapsto
  \left(e^{2\pi\sqrt{-1} u_1},\ldots, e^{2\pi\sqrt{-1}u_n}\right).$$
If $G$ is a representative of $g=(g_1,g_2,\ldots,g_n)\in T$,
and $h=\sum_{i=1}^nh_ie_i$ in $V_\Z$, then
$e^{\ll h,2\pi\sqrt{-1}G\rr}$ is equal to
$\prod_{i=1}^ng_i^{h_i}=g^h$.
As the set $\CB_n$ is not unimodular, the sets $T(\sigma)$ are not
reduced to $1$.

\begin{example} \label{exam.B2}
Let $\sigma$ be the basic set $\{e_1+e_2,e_1-e_2\}$ for $B_2$.
Then $T(\sigma)=\{(1,1),(-1,-1)\}$.
\end{example}

We now determine a set $F$ containing all sets $T(\sigma)$.

\begin{lemma} \label{lemm.B.F}
Let $\sigma$ be a basic subset of $\CB_n$.
Assume $g\in T(\sigma)$.
Then all the coordinates of $g$ are equal to $\pm 1$.
Furthermore, if $g$ is not $1$, there are at least two coordinates
of $g$ which are equal to $-1$.
\end{lemma}

\begin{proof}
We prove this by induction on $n$.
For $\CB_2$, we have seen this by direct computation.

Let $\sigma$ be a basic subset of $\CB_n$.
Assume first that $\sigma$ contains a root $e_i$.
Up to renumbering, we may assume that this root is $e_n$.
Then the basis $\sigma$ produces a basis $\sigma'$ of $\CB_{n-1}$
by putting $e_n=0$.
Let $g=(g_1,g_2,\ldots,g_n)$ in $T(\sigma)$.
We see that $g'=(g_1,g_2,\ldots,g_{n-1})$ is in $T(\sigma')$.
Thus, by induction the first $n-1$ coordinates of $g'$ are equal
to $\pm 1$.
But since $e_n$ is in $\sigma$ we get $1=g_n$.
Note that $g\neq 1$ if and only if $g'\neq 1$, hence by
induction hypothesis $g'$ has at least two coordinates not
equal to $1$.

Consider now the case where $\sigma$ does not contain any root
$e_i$.
Up to renumbering, it contains a root $e_{n-1}-e_n$ or
$e_{n-1}+e_n$.

Let us examine first the case where $\sigma$ contains the root
$\alpha=e_{n-1}-e_n$.
Let $g=(g_1,g_2,\ldots,g_{n-1},g_n)$ in $T(\sigma)$.
This implies $g_{n-1}=g_n$.
Consider the map $q$ sending $e_i$ to $e_i$ if $i<n$ and
$e_n$ to $e_{n-1}$.
Then $q$ sends $\sigma\setminus\{e_{n-1}-e_n\}$ to a basis
$\sigma'$ of $\CB_{n-1}$.
The element $g'=(g_1,g_2,\ldots,g_{n-1})$ is easily seen to
belong to $T(\sigma')$.
Indeed if $\alpha$ equals $e_i\pm e_j$ with $1\leq i<j<n$,
this is by definition.
On the other hand $q(e_i\pm e_n)=e_i\pm e_{n-1}$ and
$g_{n-1}=g_n$ imply that $g_ig_{n-1}^{\pm 1}$ coincides with
the value of $g_ig_n^{\pm 1}$.
By induction hypothesis, all coordinates  of $g'$ are equal
to $\pm 1$.
Moreover $g\neq 1$ if and only if $g'\neq 1$, so that $g$ is
of the desired form.

Finally, the same argument works if $\sigma$ contains
$\alpha=e_{n-1}+e_n$, by considering the map $q$ sending $e_i$
to $e_i$ if $i<n$, and $e_n$ to $-e_{n-1}$.
\end{proof}

\begin{definition} \label{defi.B.F}
If $I$ is a subset of $\{1,2,\ldots,n\}$ with at least two
elements, we consider the set
$F(I):=\{(g_1,g_2,\ldots,g_n)\,|\,
         g_i=-1,i\in I; g_j=1,j\notin I\}$.

We define $F\subset T$ to be the finite subset of $T$ union of
such sets $F(I)$ together with the identity $(1,1,\ldots,1)$.
\end{definition}

Let $v=\sum_{i=1}^n v_i e_i\in V$ be a vector with real
coordinates and $h=\sum_{i=1}^nh_ie_i\in V$ a vector with integral
coordinates. We will compute the normalized volume of
$\Pi_{\CB_n}(v)$ and the number of integral points in
$\Pi_{\CB_n}(h)$ using Theorem~\ref{theo.main}.
Thus we introduce the function $J_B(v)$ on $U$ defined by:

$$J_{B}(v)(u)=
  \frac{e^{\sum_{i=1}^nu_iv_i}}
       {\prod_{i=1}^nu_i
        \prod_{1\leq i<j\leq n}(u_i-u_j)
        \prod_{1\leq i<j\leq n}(u_i+u_j)}.$$
For $g=(g_1,g_2,\ldots,g_n)\in F$ and $h\in
V_\Z\cap\CC(\CB_n)$ the Kostant fraction~\eqref{equa.Fgh} is the
function on $U$ defined by:

\begin{eqnarray*}
F_{B}(g,h)(u)&=&
\frac{\prod_{i=1}^ng_i^{h_i}e^{\sum_{i=1}^n{u_ih_i}} }
     {\prod_{i=1}^n(1-g_i^{-1}e^{-u_i})
      \times
      \prod_{1\leq i<j\leq n}(1-g_i^{-1}g_je^{-(u_i-u_j)})}\\
&&\times
  \frac{1}
       {\prod_{1\leq i<j\leq
       n}(1-g_i^{-1}g_j^{-1}e^{-(u_i+u_j)})}.
\end{eqnarray*}

We have then

\begin{theorem} \label{theo.B}
Let $\gc$ be a chamber of $\CC(\CB_n)$.
\begin{itemize}
\item For any $v\in\overline{\gc}$, we have
  $$\vol_{\Z,\CB_n}(v)=\JK_{\gc}\left(J_{B}(v)\right).$$
\item For any $h\in V_\Z\cap\overline{\gc}$, the value of
  the partition function is given by:
  $$N_{\CB_n}(h)
    =\sum_{g\in F}\JK_{\gc}(F_B(g,h)).$$
\end{itemize}
\end{theorem}

As in the case of $A_n$, we will use the change of variable
$1+z_i=e^{u_i}$ to compute more easily $N_{\CB_n}(h)$. However,
let us note that this transformation does not leave the hyperplane
$u_i+u_j=0$ fixed. This hypersurface is transformed into the
hypersurface $z_i+z_j+z_iz_j=0$. So we use the expression of
$\JK_{\gc}$ as an integral over the cycle $H(\gc)$ defined in
Theorem~\ref{theo.intJK}. This cycle (its homology class) is
stable by the transformation $e^{u_i}=1+z_i$ which is close to the
identity.
Thus define the following function on $U$:

\begin{eqnarray*}
\CF_B(g,h)(z)&=&\frac
  {\prod_{i=1}^n(1+z_i)^{h_i+2n-i-1}\times\prod_{i=1}^ng_i^{h_i}}
  {\prod_{i=1}^n(1+z_i-g_i)
   \times
   \prod_{1\leq i<j\leq n}(1+z_i-g_ig_j(1+z_j))}\\
   &&\times
   \frac{1}{\prod_{1\leq i<j\leq n}(1+z_i)(1+z_j)-g_ig_j}.
\end{eqnarray*}
Performing the change of variables $e^{u_i}=1+z_i$ on the
function $F_B(g,h)(u)$ and computing the Jacobian,
Theorem~\ref{theo.main} becomes:

\begin{theorem} \label{theo.B.2}
Let $\gc$ be a chamber of $\CC(\CB_n)$.
\begin{itemize}
\item For any $v\in\overline{\gc}$, we have
  $$\vol_{\Z,\CB_n}(v)=\JK_{\gc}\left(J_{B}(v)\right).$$
\item For any $h\in V_\Z\cap\overline{\gc}$, the value of
  the partition function is given by:
  $$N_{\CB_n}(h)
    =\sum_{g\in F}\frac{1}{(2\pi\sqrt{-1})^n}\int_{H(\gc)}\CF_B(g,h)(z)dz.$$
\end{itemize}
\end{theorem}

As in the case of type $A$, in order to implement these formulae
we first have to describe the set $\CP(v,\CB_n)$
(Section~\ref{sect.B.mpns}), then  we will explain how the
integral over the cycle $H(\gc)$ is calculated similarly to an
iterated residue formula associated to these paths
(Section~\ref{sect.B.ires}), using an estimate of the order of
poles.
Finally we explain how these computations fit together to get
a global procedure for Kostant partition function for
$B_n$ (Section~\ref{sect.B.glob}).

\subsection{The search for maximal proper nested sets}
\label{sect.B.mpns}

A height function is

$$\h(v)=\sum_{i=1}^n(n+1-i)v_i$$
which takes value $1$ on all simple roots.
We will deform $\h$ later on in order to have a function taking
different values on roots.

We now proceed to describe hyperplanes for $\CB_n$. If
$P=[P^+,P^-]$ are two disjoints subsets of $\{1,2,\ldots,n\}$, we
denote by $\ll u_P,v\rr$ the linear form $\sum_{i\in
P^+}v_i-\sum_{j\in P^-}v_j$. Consider the hyperplane

$$H_P=\{v\in V,\ll u_P,v\rr=0\}$$
in $V$.
It is equal to the hyperplane determined by the reverse
list $[P^-,P^+]$.
Thus to each set $P=\{P^+,P^-\}$ of two disjoint sets $P^+$, $P^-$
such that at least one is non empty, we associate a hyperplane $H_P$.

We denote by $Z$ the complement of $P^+\cup P^-$ in
$\{1,2,\ldots,n\}$ and by $\CB(Z)$ the subset of $\CB_n$
defined by

$$\CB(Z)=
      \{e_i\,|\,i\in Z\}
  \cup\{e_i\pm e_j\,|\,1\leq i<j\leq n\,;\,i,j\in Z\}.$$
This is the positive root system $B_{|Z|}$,
with the positivity induced by the lexicographic order.

Let $\CK(P^+,P^-)$ be the subset of $\CB_n$ defined by

\begin{eqnarray*}
&    &\{e_i-e_j\,|\,1\leq i<j\leq n\,;\,i,j\in P^+\}\\
&\cup&\{e_i+e_k\,|\,i\in P^+,k\in P^-\}\\
&\cup&\{e_k-e_\ell\,|\,1\leq k<\ell\leq n;\,k,\ell\in P^-\}.
\end{eqnarray*}
Note that by defining $f_i=e_i$ if $i\in P^+$ and
$f_k=-e_{|P^-|-k+1}$ if $k\in P^-$, the set $\CK(P^+,P^-)$
coincides with

\begin{eqnarray*}
&    &\{f_i-f_j\,|\,1\leq i<j\leq n;\,i,j\in P^+\}\\
&\cup&\{f_i-f_k|\,\,i\in P^+,k\in P^-\}\\
&\cup&\{f_k-f_\ell\,|\,1\leq k<\ell\leq n;\,k,\ell\in P^-\}.
\end{eqnarray*}
Thus the set $\CK(P^+,P^-)$ is a positive root system of
type $A_{|P^+|+|P^-|-1}$.
However the positivity is induced by the lexicographic order
on $P^+$ and the reverse lexicographic order on $P^-$.
Observe also that $H_P$ is the vector space spanned by
$\CK(P^+,P^-)\cup\CB(Z)$.

\begin{lemma} \label{lemm.B}
\begin{itemize}
\item The hyperplane $H_P$ is a $\CB_n$-admissible hyperplane.
\item The set $\CB_n\cap H_P$ is the union of $\CB(Z)$ and
  $\CK(P^+,P^-)$.
\item Every $\CB_n$-admissible hyperplane is of this form.
\end{itemize}
\end{lemma}

\begin{proof}
The first two assumptions are easy to see.
We prove the third assumption by induction on $n$, the case
$n=2$ being trivial.
Let $H$ be a $\CB_n$-admissible hyperplane
and let $\alpha$ be a root in $H$.
There are $3$ possibilities for $\alpha$: up to renumbering roots,
we can consider the cases $\alpha=e_n$, $\alpha=e_{n-1}-e_n$
and $\alpha=e_{n-1}+e_n$.

In the first case, the map $q$ sending $e_i$ to $e_i$ if $i<n$ and
$e_n$ to $0$ maps the set $\CB_n\setminus\{\alpha\}$ to
$\CB_{n-1}$.
The space $H/\R\alpha$ becomes a $\CB_{n-1}$-admissible
hyperplane.
It is thus determined by $P'=[P^{'+}, P^{'-}]$, where
$P^{'+}$ and $P^{'-}$ are two disjoint sets
contained in $\{1,2,\ldots, n-1\}$.
Then the hyperplane $H$ is equal to the hyperplane determined
by $[P'^+,P'^-]$.

In the second case, the map $q$ sending $e_i$ to $e_i$ if $i<n$
and $e_n$ to $e_{n-1}$ sends the set $\CB_n\setminus\{\alpha\}$
to $\CB_{n-1}$.
The space $H/\R\alpha$ becomes a $\CB_{n-1}$-admissible
hyperplane.
It is thus determined by $P'=[P^{'+} ,P^{'-}]$.
If neither $P^{'+}$ nor $P^{'-}$ contain $n-1$, the
hyperplane $H$ is equal to the hyperplane determined by
$[P^{'+},P^{'-}]$.
Otherwise assume that for example $P^{'+}$ contains $n-1$.
Then the hyperplane $H$ is equal to the hyperplane
determined by $[P^+,P^-]$, where $P^+=P^{'+}\cup\{n\}$
and $P^-=P^{'-}$.

In the third case, the map $q$ sending $e_i$ to $e_i$ if $i<n$
and $e_n$ to $-e_{n-1}$ sends the set
$\CB_n\setminus\{\alpha\}$ to $\CB_{n-1}$.
The space $H/\R\alpha$ becomes a $\CB_{n-1}$-admissible
hyperplane.
It is thus determined by $P'=[P^{'+}, P^{'-}]$.
If neither $P^{'+}$ nor $P^{'-}$ contains $n-1$, the hyperplane
$H$ is equal to the hyperplane determined by $[P^{'+},P^{'-}]$.
Assume that $P^{'+}$ contains $n-1$.
Then the hyperplane $H$ is equal to the the hyperplane
determined by $[P^+,P^-]$, where $P^+=P^{'+}$ and
$P^-=P^{'-}\cup\{n\}$.
\end{proof}

We now give a detailed description of our algorithm computing
maximal nested sets.
We describe a vector as an array $v=[v_1,v_2,\ldots,v_n]$. To
check if $v$ is in the cone $\CC(\CB_n)$, we need to verify if
$\sum_{j=1}^iv_j\geq 0$ for $1\leq i\leq n$. This is done by the
procedure ${\tt CheckBvector}$, which returns the answer true or
false.

For the system $\CB_n$ the highest root $\theta^B(n)$ is equal to

$$\theta^B(n)=[1,1,0,0,0,\ldots,0].$$

We recall here that $P$ is divided in two sets $P^+\cup P^-$, one
of them being non empty. The first task is to list the hyperplanes
in $\Hyp(v,\CB_n)$. This set of hyperplanes is obtained by the
command line ${\tt AllPossibleBwalls(v)}$. The input of this
procedure is the vector $v$. The output is a set of elements
$P=\{P^+,P^-\}$, where $P^+=[i_1,i_2,\ldots i_p]$ and
$P^-=[j_1,j_2,\ldots, j_q]$ are two ordered disjoint lists made
from indices taken in $\{1,\ldots,n\}$, with at least one of $P^+$
or $P^-$ being non empty. Let $\ll u_P,v\rr=\sum_{i\in
P^+}v_i-\sum_{j\in P^-}v_j$ be the normal vector to $H_P$. Then as
stated in Lemma~\ref{lemm.cond} we need to test if $\ll
u_P,\theta^B(n)\rr$ is not zero and if $\ll u_P,v\rr\times\ll
u_P,\theta^B(n)\rr$ is non negative.

We then construct the vector
$$\proj_H(v)=v
-\frac{\ll u_P,v\rr}
      {\ll u_P,\theta^B(n)\rr}\theta^B(n).$$
This vector is represented as $\{{[v_1,P^+],[v_2,P^-]},[w,Z]\}$.
The sum of coordinates of $v_1$ is equal to the sum of the
coordinates of $v_2$.
Now $Z$ is the ordered list $[k_1,k_2,\ldots,k_\ell]$ of
complementary indices to $P^+\cup P^-$ and

$$w=[\proj_H(v)[k_1],\ldots,\proj_H(v)[k_\ell]].$$

Note that the equations of the cone $\CC(\CK(P^+,P^-))$ can
be given in the convenient form
$v_1\oplus v_2\in\CC(\CK(P^+,P^-))$
if and only if ${\tt CheckBvector(v_1)}$ and
${\tt CheckBvector(v_2)}$ are true.
Equations of the cone $\CC(\CB(Z))$ are given in the form
$w\in\CC(\CB(Z))$ if and only if ${\tt CheckBvector(w)}$ is true.

Thus the condition that $H$ is in $\Hyp(v,\CB_n)$ is equivalent to
the series of conditions:

\begin{eqnarray*}
\ll u_P,\theta^B(n)\rr&\neq&0,\\
\ll u_P,v\rr\times\ll u_P,\theta^B(n)\rr&\geq&0,\\
{\tt CheckBvector(v_1)}&=&{\tt true},\\
{\tt CheckBvector(v_2)}&=&{\tt true},\\
{\tt CheckBvector(w)}&=&{\tt true}.
\end{eqnarray*}
Those five conditions are checked by the command line
${\tt CheckBwall(v,H)}$, that gives an answer true or false.

\begin{remark} \label{rema.bwall}
We can first construct all disjoint subsets $P^+$,
$P^-$ of $\{1,2,\ldots,n\}$ and test these five conditions
successively on all of them.
However it is highly desirable to throw away \emph{a priori} a
great number of these partitions by noticing the following
restrictive conditions on the possible lists to be considered.

Let $\{P^+,P^-\}=\{[i_1,i_2,\ldots i_p],[j_1,j_2,\ldots,j_q]\}$
be a set of two disjoint subset of $\{1,2,\ldots n\}$ represented
as lists with strictly increasing indices.
Let $Z=[k_1,k_2\ldots,k_\ell]$ be the list of complementary
indices to $P^+\cup P^-$ in $\{1,\ldots,n\}$.
The following linear forms are positive on the cone
$\CC(\CK(P^+,P^-))$ generated by $\CK(P^+)$ and $\CK(P^-)$:

\begin{eqnarray*}
v_{i_1}+v_{i_2}+\cdots+v_{i_s}&\geq&0
  \quad\quad\mbox{for all $1\leq s\leq p$,}\\
v_{j_1}+v_{j_2}+\cdots+v_{j_t}&\geq&0
  \quad\quad\mbox{for all $1\leq t\leq q$,}\\
v_{k_1}+v_{k_2}+\cdots+v_{k_s}&\geq&0
  \quad\quad\mbox{for all $1\leq s\leq\ell$.}\\
\end{eqnarray*}

Note that all the above linear forms take positive values on
$\theta^B(n)$.
We employ Lemma~\ref{lemm.B}.
Thus if $v[i_1]<0$, the index $i_1$ cannot start the list
$P^+$ of an element $\{P^+,P^-\}$ in ${\tt AllPossibleBwalls(v)}$
and we reject all such $\{P^+,P^-\}$.

Similarly assume that we have constructed a list of indices
$[i_1,i_2]$ satisfying conditions $v[i_1]\geq 0$ and
$v[i_1]+v[i_2]\geq 0$.
Then if $v[i_1]+v[i_2]+v[i_3]<0$, a list starting with
$[i_1,i_2,i_3]$ cannot be the first three indices of the
component $P^+$ of an element $\{P^+,P^-\}$ in the set
${\tt AllPossibleBwalls(v)}$ and we skip it right away.
\end{remark}

This achieves the description of the procedure
${\tt AllPossibleBwalls}$.
We now have to perform the next step of our algorithm.
As for type $A$ we build MNSs iteratively.
At each step we get a set of partial MNSs, to which we will
apply recursively our algorithm.
Note that after Lemma~\ref{lemm.B} the intersection of
a $\CB_n$-admissible hyperplane $H_P$ with $\CB_n$ is the
union of a system of type $A$ and a system of type $B$.

The part of the MNS coming from the subsystem of type $A$
is computed with the procedure ${\tt AddAnests}$.
It performs a reordering of the result of a call to the
procedure ${\tt MNSs}$ described in Section~\ref{sect.A.mpns}.

The part of the MNS coming from the subsystem of type $B$
is computed with the procedure ${\tt Bsplits}$, calling the
previously described procedure ${\tt AllPossibleBwalls}$.

Procedures ${\tt AddAnests}$ and ${\tt Bsplits}$ are enclosed
in ${\tt MoreNSs}$, thus giving a new iteration of the process.
After regularization of the result we hence get a procedure
named ${\tt AllNewNSs}$, performing a new step in the building
of MNSs.

Finally the procedure ${\tt B\_MNSs}$, computing MNSs for a
given vector $v$ for type $B$, is the following.
First, we use a procedure named ${\tt B\_NS1}$ to calculate
the first seed of all MNSs.
After, repeated calls to the procedure ${\tt AllNewNSs}$
build the desired MNSs.

\subsection{Residues associated to maximal proper nested sets}
\label{sect.B.ires}

A proper maximal nested set $M$ gives rise to an ordered basis
$\alpha_i$, and a cycle  $H(M)$.
We need to compute

$$\int_{H(M)}\CF_B(g,h)(z)dz$$
where
$$H(M):=\{z, |\ll\alpha_i,z\rr|=\epsilon_i\}.$$

The function $z\mapsto\CF_B(g,h)(z)$ is deduced from the function
$F_{B}(g,h)(u)$ in the space $\widehat R_{\CA}$ by the change of
variable $e^{u_i}=1+z_i$.
Thus its denominator is a product of factors, either of the
form $z_i$ corresponding to the root $u_i$, or of the
form $z_i-z_j$ corresponding to the root $u_i-u_j$ or
$z_i+z_j+z_iz_j$ corresponding to the root $u_i+u_j$.
We denote by $u(z)$ the point with  coordinates $u_i$ satisfying
$e^{u_i}=1+z_i$.

We start integrating our function $\CF_B(g,h)(z)$ over the smaller
circle $|\ll\alpha_1,z\rr|=\epsilon_1$ keeping the other variables
fixed. By our condition on the cycle, the function we integrate
has poles on the domain $|\ll\alpha_1,z\rr|\leq\epsilon_1$ only
when $\alpha_1(u(z))=0$. If $\alpha_1(u(z))=u_i-u_j$ or
$\alpha_1(u(z))=u_i$, the poles are obtained for $z_i=z_j$ or
$z_i=0$. If $\alpha_1(u(z))=u_i+u_j$, the pole on the domain
$|\ll\alpha_1,z\rr|\leq\epsilon_1$ is obtained for
$z_i=-z_j/(1+z_j)$. Thus we compute the integral over the circle
by the residue theorem in one variable, and proceed. From the
general theory, the poles of the function we obtain, replacing
$z_i$ by one of the values above are again of the same form with
respect to the remaining variables, as is easily checked.

As in case $A_{n-1}$, for a root $\alpha=u_i$ (resp.
$\alpha=u_i\pm u_j$) we can replace after taking the residue at
$\alpha=0$ the variable $z_i$ by $0$ (resp. by $\mp z_j$) in all
other roots.
Thus we get rid of the variable $z_i$.
The procedure ${\tt FormalPathB}$ produces the ordered
path resulting from all these substitutions.

In the case of type $B$ we compute the residue by directly checking the
order of the pole at $\alpha=0$, and then using differentiation.
The program works in the same way with parameters. The function
obtained is locally polynomial with polynomial coefficients
depending of the parity of the integers $h_i$.

\subsection{The procedure ${\tt MNS\_KostantB}$}
\label{sect.B.glob}

We finish the section dedicated to $B_n$ by giving the global
outline of the procedure ${\tt MNS\_KostantB(v)}$ computing the
Kostant partition number of a vector $v$ lying in the root lattice
of $\CB_n$. We begin by slightly deforming $v$ so that it lies on
no wall, by setting ${\tt v':=DefVectorB(v,n)}$. We then compute
all MNSs for $v'$ with the call ${\tt B\_MNSs(v')}$
(Section~\ref{sect.B.mpns}). For every MNS $M$, we extract the
list $R$ of highest roots of its irreducible components by setting
${\tt R:=BthetaMNS(M)}$. We sort these roots by their height with
the command line ${\tt R':=BorderThetas(R,n)}$.
We then transform the list of roots $R'$ into a path $p$ by setting
${\tt p:=FormalPathB(R')}$.

Now remark that our procedures are designed to take residues along
positive roots, using the fact that $\res_{-\alpha}=-\res_{\alpha}$
for any root $\alpha$.
The sign that appears (more precisely $-1$ to the power the
number of negative roots in the path $p$) is computed with the
procedure ${\tt PathSign(p,n)}$.

Then for every $g$ in $F$ we do the following.
The iterated residue along the path $p$ and for $g$ is obtained by the
command line ${\tt OneIteratedBresidue(p,g,v,n)}$.
Let us briefly describe its implementation.
We first compute the Kostant fraction (second item of
Definition~\ref{defi.func}, procedure ${\tt KostantFunctionB}$).
Then for every root of the path we apply the procedure
${\tt ComputeOneResidue}$ (Section~\ref{sect.B.ires}) and update
the order of the pole with a procedure named ${\tt OrderPoleB}$.

Finally summing all products
${\tt PathSign(p,n)}\times{\tt OneIteratedBresidue(p,g,v,n)}$
over the sets of $g$'s and of $M$'s, we get the
desired result.

\begin{remark} \label{rema.comp}
Let us fix a list $R'=[\alpha_1,\ldots,\alpha_n]$ of ordered
roots coming from a MNS, and an element $g$.
We say that $R'$ and $g$ are \emph{compatible} if the following
condition is satisfied.
If indices of monomial(s) of $\alpha_k$ have not yet occured
among indices of roots $\alpha_\ell$ with $\ell<k$, then
$g$ must satisfy $g^{\alpha_k}=1$
(that is $g_ig_j^{\pm 1}=1$ if $\alpha_k=e_i\pm e_j$ and
$g_i=1$ if $\alpha_k=e_i)$.
Note that the iterated residue for $g$ and for the path $p$
associated to $R'$ is zero if $g$ and $R'$ are not compatible.
Hence summing only over $g$'s that are compatible with a given
list $R'$ saves useless computations.
The check of compatibility is performed by the procedure
${\tt ListAndGAreCompatible(R',g,n)}$.
\end{remark}

\section{The type $C_n$}
\label{sect.C}

Consider a vector space $V$ with basis $e_1$, $e_2$, \ldots,
$e_n$. We choose on $V$ the standard Lebesgue measure $dh$.
Let
$$\CC_n=\{2e_i\,|\,1\leq i\leq n\}
        \cup\{e_i-e_j\,|\,1\leq i<j\leq n\}
        \cup\{e_i+e_j\,|\,1\leq i<j\leq n\}.$$
Then $\CC_n$ is a positive roots system of type $C_n$,
and generates $V$.
The number of elements in $\CC_n$ is $N=n^2$.
Note that elements of $\CC_n$ and $\CB_n$ are proportional, so they
determine the same hyperplane arrangement and the same chambers.

Let $L$ be the lattice defined by $\Z e_1\oplus\Z
e_2\oplus\cdots\oplus\Z e_n$. We remark that the lattice $V_\Z$
generated by $\CC_n$ is the sublattice of index $2$ in $L$
consisting of all elements $v=[v_1,v_2,\ldots,v_n]$ with integral
coordinates and such that the sum $\sum_{i=1}^n v_i$ is an even
integer. A $\Z$-basis of $V_\Z$ is, for example,
$$\Z(e_1-e_n)\oplus\Z(e_2-e_n)\oplus\cdots\oplus
  \Z(e_{n-1}-e_n)\oplus\Z(2e_n) \ ,$$
so $\vol \left( V/V_\Z \right) = 2$.

The dual lattice $U_\Z$ is the lattice of vectors
$\gamma=(\gamma_1,\gamma_2,\ldots,\gamma_n)$ such that $\gamma_i$
are half integers and such that $\gamma_i+\gamma_j$ is an integer
for all $i$, $j$.
The set $U_\Z/\Z e_1\oplus\cdots\oplus\Z e_n$ is of cardinality
$2$ with representative elements $(0,0,\ldots,0,0)$ and
$(1/2,\ldots,1/2)$.

As before, we identify the torus
$\tilde{T}=U/(\Z e_1\oplus\cdots\oplus\Z e_n)=\R^n/\Z^n$
with $(S^1)^n$ by

$$(u_1,u_2,\ldots, u_n)\mapsto
  \left(e^{2\pi\sqrt{-1} u_1},\ldots, e^{2\pi\sqrt{-1}u_n}\right).$$
 Then  $$T=\tilde{T}/\{\pm 1\}=U/U_\Z.$$
Let $G$ be a
representative of $g=(g_1,g_2,\ldots,g_n)\in \tilde{T}$ and
$h=\sum_{i=1}^nh_ie_i$ in $V_\Z$. Then $e^{\ll
h,2\pi\sqrt{-1}G\rr}$ is equal to $\prod_{i=1}^ng_i^{h_i}=g^h$.
This function  is well defined on $T=\tilde{T}/\{\pm 1\}$ since
$\sum_{i=1}^nh_i$ is even.

For $\sigma$ a basic subset of $\CC_n$, define

$$\tilde{T}(\sigma)=\left\{g\in\tilde{T}\,\Big|\,
e^{\ll\alpha,2\pi\sqrt{-1}G\rr}=1\,\,
  \mbox{for all}\,\alpha\in\sigma\right\}.$$
As the set $\CC_n$ is  not unimodular, sets $\tilde{T}(\sigma)$
are not reduced to $1$.

\begin{lemma} \label{lemm.C}
Let $\sigma$ be a basic subset of $\CC_n$.
Then $\tilde{T}(\sigma)\subset\{\pm 1\}^n$.
\end{lemma}

\begin{proof}
We prove by induction on $n$ that if $\sigma$ is basic then
the condition $g=(g_1,\ldots,g_n)\in \tilde{T}(\sigma)$ forces
$g_i^2=1$ ($1\leq i\leq n$).
In other words $g^\alpha=1$ for all long roots $\alpha$.
If so then $g_i=\pm 1$ for all $i$.
The base of the induction, that is $\CC_2$, is straightforward
and we omit it.
We thus proceed considering various possibilities for our $\sigma$.

If there exists a long root in $\sigma$ we may assume that this
long root is $2e_n$.
We embed the system $\CC_{n-1}$ in $\CC_n$ via the first
$(n-1)$ coordinates.
Then the basis $\sigma$ of $\CC_n$ produces a basis $\sigma'$
of $\CC_{n-1}$ consisting of roots
$\{e_i\pm e_j\in\sigma\,|\,1\leq i<j\leq n-1\}$, of roots
$\{2e_i\in\sigma\,|\,1\leq i\leq n-1\}$,
and of roots $\{2e_i\,|\,e_i\pm e_n\in\sigma;i\neq n\}$.
It is easy to see that the elements $(g_1,g_2,\ldots,g_{n-1})$
are in $\tilde{T}(\sigma')$.
Indeed $g_i^2=1$ if $e_i\pm e_n\in\sigma$ as $g_ig_n^{\pm 1}=1$
and $g_n^2=1$; and similarly $g_i^2=1$ if $2e_i\in\sigma$.
Thus by induction we obtain $g_i^2=1$ for every $i$.

Now assume that there is no long root in $\sigma$.
We may assume that there is a root of the form
$e_{n-1}-e_n$ or $e_{n-1}+e_n$.

In the first case, consider the basis $\sigma'$ of $\CC_{n-1}$
consisting
of the roots $\{e_i\pm e_j\in\sigma\,|\,1\leq i<j\leq n-1\}$ and
of the roots $\{e_i\pm e_{n-1}\,|\,e_i\pm e_n\in\sigma\}$.
It is easy to see that the elements $(g_1,g_2,\ldots,g_{n-1})$
are in $\tilde{T}(\sigma')$.
Indeed, for example, $g_ig_{n-1}^{\pm 1}=1$ if
$e_i\pm e_{n-1}\in\sigma'$, as $g_ig_n^{\pm 1}=1$ and
$g_n=g_{n-1}$.
Thus by the induction hypothesis we obtain $g_i^2=1$ for all
$i\neq n$.
Since $g_n=g_{n-1}$, we also obtain $g_n^2=1$.

The second case is similar.
\end{proof}

Let  $v=\sum_{i=1}^n v_i e_i\in V$ be a vector with real
coordinates and $h=\sum_{i=1}^nh_ie_i\in V$ a vector with integral
coordinates and such that $\sum_{i=1}^nh_i$ is even. We will
compute the normalized volume of $\Pi_{\CC_n}(v)$ and the number
of integral points in $\Pi_{\CC_n}(h)$ using
Theorem~\ref{theo.main}. We will use the JK residue with respect
to the  measure $dh$ associated to the basis $e_1,e_2,\ldots,
e_n$. However, the normalized volume $\vol_{\Z,\CC_n}(h)$  is
computed for the measure determined by the lattice spanned by
$\CC_n$ which is of index $2$ in $\oplus_{i=1}^n \Z e_i$.

We introduce the function $J_C(v)$ on $U$ defined by:

$$J_{C}(v)(u)=
  \frac{e^{\sum_{i=1}^nu_iv_i}}
       {\prod_{i=1}^n2u_i
        \prod_{1\leq i<j\leq n}(u_i-u_j)
        \prod_{1\leq i<j\leq n}(u_i+u_j)}.$$

For $g=(g_1,g_2,\ldots,g_n) \in\{\pm 1\}^n$ the Kostant
fraction~\eqref{equa.Fgh} is the function on $U$ defined by:

\begin{eqnarray*}
F_C(g,h)(u)&=& \frac{\prod_{i=1}^ng_i^{h_i}e^{\sum_{i=1}^nu_ih_i}}
     {\prod_{i=1}^n(1-e^{-2u_i})
      \times
      \prod_{1\leq i<j\leq n}(1-g_i^{-1}g_je^{-(u_i-u_j)})}\\
&&\times
  \frac{1}
       {\prod_{1\leq i<j\leq n}(1-g_i^{-1}g_j^{-1}e^{-(u_i+u_j)})}.
\end{eqnarray*}

\begin{theorem} \label{theo.C}
Let $\gc$ be a chamber of $\CC(\CC_n)$.
\begin{itemize}
\item For any $v\in\overline{\gc}$, we have
  $$\vol_{\Z,\CC_n}(v)=2\,\JK_{\gc}\left(J_{C}(v)\right).$$
\item For any vector $h\in V_\Z\cap\overline{\gc}$ with integral
  coordinates such that $\sum_{i=1}^nh_i$ is even, the value of
  the partition function is given by:
  $$N_{\CC_n}(h)
    =\sum_{g\in \{\pm 1\}^n}\JK_{\gc}(F_C(g,h)).$$
\end{itemize}
\end{theorem}

In the second formula, there should be a multiplication by a
factor  $2$ as the volume of the fundamental domain of the lattice
spanned by $\CC_n$ is $2$. However, we should  sum only on
$T=\tilde T/\{\pm 1\}$. Thus the two factors of $2$ compensate each other. In
fact,  we will indeed  sum over  $T$ represented as $\{\pm
1\}^{n-1}\times\{1\}$ and multiply the result by the constant $2$.

As in the case of $B_n$, we will use the change of variable
$1+z_i=e^{u_i}$ to compute more easily the formula for
$N_{\CC_n}(h)$. As explained in the case of $B_n$ we need to use
the integral formulation of the Jeffrey-Kirwan residue.
Thus define

\begin{eqnarray*}
\CF_C(g,h)(z)&=&\frac
  {\prod_{i=1}^n(1+z_i)^{h_i+2n-i}\times\prod_{i=1}^ng_i^{h_i}}
  {\prod_{i=1}^n((1+z_i)^2-1)
   \times
   \prod_{1\leq i<j\leq n}(1+z_i-g_ig_j(1+z_j))}\\
   &&\times
   \frac{1}{\prod_{1\leq i<j\leq n}(1+z_i)(1+z_j)-g_ig_j}.
\end{eqnarray*}
Performing the change of variables $e^{u_i}=1+z_i$ on
the function $F_C(g,h)(u)$ and computing the Jacobian,
Theorem~\ref{theo.main} becomes:

\begin{theorem} \label{theo.C.bis}
Let $\gc$ be a chamber of $\CC(\CC_n)$.
\begin{itemize}
\item For any $v\in\overline{\gc}$, we have
  $$\vol_{\Z,\CC_n}(v)=2\,\JK_\gc\left(J_{C}(v)\right).$$
\item For any vector $h\in V_\Z\cap\overline{\gc}$ with integral
  coordinates $h_i$ with $\sum_{i=1}^n h_i$ even, the value of
  the partition function is given by:
  $$N_{\CC_n}(h)
    =\sum_{g\in \{\pm 1\}^n} \frac{1}{(2\pi\sqrt{-1})^n}\int_{H(\gc)}\CF_C(g,h)(z)dz.$$
\end{itemize}
\end{theorem}

Similarly we will  sum over  $T$ represented as $\{\pm
1\}^{n-1}\times\{1\}$ and multiply the result by the constant $2$.

The cycle $H(\gc)$ associated to a chamber $\gc$  containing a
regular element $v=\sum_{i=1}^n v_i e_i$ is the same cycle that we
computed in the preceding section for $B_n$ . Hence we can reuse
most of procedures from the type $B_n$. Paths are the same, and
the residue calculations are the same. More precisely, the only
two changes are in the computation of the set $G$ (procedure ${\tt
GC(n)}$) and in the computation of the Kostant function (procedure
${\tt UCKostant}$). This terminates the case of $C_n$.

\section{The type $D_n$}
\label{sect.D}

\subsection{The formulae to be implemented}
\label{sect.D.form}

Consider a vector space $V$ with basis $e_1$, $e_2$, \ldots,
$e_n$. We choose the standard Lebesgue measure $dh$.
Let

$$\CD_n=\{e_i-e_j\,|\,1\leq i<j\leq n\}
        \cup\{e_i+e_j\,|\,1\leq i<j\leq n\}.$$
Then $\CD_n$ is a positive roots system of type $D_n$, and
generates $V$.
The number of elements in $\CD_n$ is $N=n^2-n$.

We remark that the lattice $V_\Z$ generated by roots  of $\CD_n$
is the same lattice as the one generated by the roots of
$\CC_n$. It is of  index $2$ in $L:=\Z e_1\oplus\Z
e_2\oplus\cdots\oplus\Z e_n$  and consists of elements
$v=[v_1,v_2,\ldots,v_n]$ with integral coordinates such that
the sum $\sum_{i=1}^n v_i$ is an even integer.
 The group $T=U/U_\Z$ is thus the quotient of $\tilde T=U/\Z
e_1\oplus\cdots\oplus\Z e_n$, obtained by identifying $g$ and
$-g$, that is $T=\tilde{T}/\{\pm 1\}$. As in Section~\ref{sect.C},
we identify the torus $\tilde{T}=U/(\Z e_1\oplus\cdots\oplus\Z
e_n)=\R^n/\Z^n$ to $(S^1)^n$ by
$$(u_1,u_2,\ldots, u_n)\mapsto
  \left(e^{2\pi\sqrt{-1} u_1},\ldots, e^{2\pi\sqrt{-1}u_n}\right).$$
Consider the set $F=\{\pm 1\}^n\subset(S^1)^n$. For $\sigma$ a
basic subset of $\CD_n$, define

$$\tilde{T}(\sigma)=\left\{g\in\tilde{T}\,\Big|\,
e^{\ll\alpha,2\pi\sqrt{-1}G\rr}=1\,\,
  \mbox{for all}\,\alpha\in\sigma\right\}.$$

\begin{lemma} \label{lemm.D.F}
Let $\sigma$ be a basic subset of $\CD_n$. Then
$\tilde{T}(\sigma)$ is contained in $F$.
\end{lemma}

\begin{proof}
Basic subsets of $\CD_n$ are basic subsets of $\CC_n$ so that we
can choose the same set $F=\{\pm 1\}^n$.
\end{proof}

Let $v=\sum_{i=1}^n v_i e_i\in V$ be a vector with real
coordinates and $h=\sum_{i=1}^nh_ie_i\in V$ a vector with integral
coordinates and such that $\sum_{i=1}^nh_i$ is even. We will
compute the normalized volume of $\Pi_{\CD_n}(v)$ and the number
of integral points in $\Pi_{\CD_n}(h)$ using
Theorem~\ref{theo.main}.

Thus we introduce the function $J_D(v)$ on $U$ defined by:

$$J_{D}(v)(u)=
  \frac{e^{\sum_{i=1}^nu_iv_i}}
       {\prod_{1\leq i<j\leq n}(u_i-u_j)
        \prod_{1\leq i<j\leq n}(u_i+u_j)}.$$

For $g=(g_1,g_2,\ldots,g_n) \in\{\pm 1\}^n$ the Kostant
fraction~\eqref{equa.Fgh} is the function on $U$ defined by:
\begin{eqnarray*}
F_D(g,h)(u)&=& \frac{\prod_{i=1}^ng_i^{h_i}\times e^{\sum_{i=1}^nu_ih_i}}
{\prod_{1\leq i<j\leq n}(1-g_i^{-1}g_je^{-(u_i-u_j)})}\\
&&\times
  \frac{1}
       {\prod_{1\leq i<j\leq n}(1-g_i^{-1}g_j^{-1}e^{-(u_i+u_j)})}.
\end{eqnarray*}

We have then

\begin{theorem} \label{theo.D}
Let $\gc$ be a chamber of $\CC(\CD_n)$.
\begin{itemize}
\item For any $v\in\overline{\gc}$, we have
  $$\vol_{\Z,\CD_n}(v)=2\,\JK_{\gc}\left(J_{D}(v)\right).$$
\item For any vector $h\in V_\Z\cap\overline{\gc}$ with integral
  coordinates such that $\sum_{i=1}^nh_i$ is even, the value of
  the partition function is given by:
  $$N_{\CD_n}(h)
    =\sum_{g\in \{\pm 1\}^n}\JK_{\gc}(F_D(g,h)).$$
\end{itemize}
\end{theorem}

We use the change of variable $1+z_i=e^{u_i}$ to compute more
easily the formula for $N_{\CD_n}(h)$ and thus introduce
integration over a cycle.
Thus define

\begin{eqnarray*}
\CF_D(g,h)(z)&=&\frac
  {\prod_{i=1}^n(1+z_i)^{h_i+2n-i-2}
   \times\prod_{i=1}^ng_i^{h_i}}
  {\prod_{1\leq i<j\leq n}(1+z_i-g_ig_j(1+z_j))}\\
   &&\times
   \frac{1}{\prod_{1\leq i<j\leq n}(1+z_i)(1+z_j)-g_ig_j}.
\end{eqnarray*}
After performing the change of variables $e^{u_i}=1+z_i$
on the function $F_D(g,h)(u)$ and after computing the Jacobian,
Theorem~\ref{theo.main} becomes:

\begin{theorem} \label{theo.D.bis}
Let $\gc$ be a chamber of $\CC(\CD_n)$.
\begin{itemize}
\item For any $v\in\overline{\gc}$, we have
  $$\vol_{\Z,\CD_n}(v)=2\,\JK_{\gc}\left(J_{D}(v)\right).$$
\item For any vector $h\in V_\Z\cap\overline{\gc}$ with integral
  coordinates $h_i$ such that $\sum_{i=1}^nh_i$
  is even, the value of the partition function is given by:
  $$N(\CD_n,h)
    =\sum_{g\in \{\pm 1\}^n}
     \frac{1}{(2\pi\sqrt{-1})^n}\int_{H(\gc)}\CF_D(g,h)(z)dz.$$
\end{itemize}
\end{theorem}

As for types $A$ and $B$, in order to implement these formulae we
first have to describe the set $\CP(v,\CD_n)$
(Section~\ref{sect.D.mpns}).
We finish to explain the implementation of case $D$ in
Section~\ref{sect.D.glob}, using the fact that types $B$ and $D$
are similar.

\subsection{The search for maximal proper nested sets}
\label{sect.D.mpns}

A height function is
$$\h(v)=\sum_{i=1}^n(n-i)v_i$$
which takes value $1$ on all simple roots.
We will deform it later on in order to have a function taking
different values on roots.

We now proceed to describe hyperplanes for $D_n$.
If $P=[P^+,P^-]$ are two disjoints subsets of $\{1,2,\ldots,n\}$,
we denote by $\ll u_P,v\rr$ the linear form
$\sum_{i\in P^+}v_i-\sum_{j\in P^-}v_j$.
Consider the hyperplane in $V$ defined by

$$H_P=\{v\in V,\ll u_P,v\rr=0\}$$
and remark that it is equal to the hyperplane determined by the
reverse list $[P^-,P^+]$.
Thus to each set $P=\{P^+,P^-\}$ of two disjoint sets $P^+$, $P^-$
such that at least one is non empty, is associated a hyperplane
$H_P$.

We denote by $Z$ the complement of $P^+\cup P^-$ in
$\{1,2,\ldots,n\}$ and by $\CD(Z)$ the subset of $\CD_n$
defined by

$$\CD(Z)=\{e_i\pm e_j\,|\,1\leq i<j\leq n;\,i,j\in Z\}.$$
This is the positive roots system of type $D_{|Z|}$,
with the positivity induced by the lexicographic order.

Let $\CK(P^+,P^-)$ be the subset of $\CD_n$ defined by

\begin{eqnarray*}
&    &\{e_i-e_j\,|\,1\leq i<j\leq n;\,i,j\in P^+\}\\
&\cup&\{e_i+e_k\,|\,i\in P^+,k\in P^-\}\\
&\cup&\{e_k-e_\ell\,|\,1\leq k<\ell\leq n;\,k,\ell\in P^-\}.
\end{eqnarray*}
As we observed in Section~\ref{sect.B.mpns} for $B_n$, by
defining $f_i=e_i$ if $i\in P^+$ and $f_k=-e_{|P|-k+1}$ if
$k\in P^-$, the set $\CK(P^+,P^-)$ is a positive roots system
of type $A_{|P^+|+|P^-|-1}$.
Here the positivity is induced by the lexicographic order on
$P^+$ and the reverse lexicographic order on $P^-$.

Observe also that $H_P$ is the vector space spanned by
$\CK(P^+,P^-)\cup\CD(Z)$.

\begin{lemma} \label{lemm.D}
\begin{itemize}
\item The hyperplane $H_P$ is a $\CD_n$-admissible hyperplane.
\item The set $\CD_n\cap H_P$ is the union of $\CD(Z)$ and
  $\CK(P^+,P^-)$.
\item Every $\CD_n$-admissible hyperplane is of this form.
\end{itemize}
\end{lemma}

\begin{proof}
The first two assumptions are easy to see. Now as $\CD_n$ is
contained in $\CB_n$, a $\CD_n$-admissible hyperplane is $\CB_n$
admissible, so is of this form.
\end{proof}

\subsection{The procedure ${\tt MNS\_KostantD}$}
\label{sect.D.glob}

Most of procedures from type $B_n$ are kept unchanged.
More precisely, the iterated residue calculation, the estimate of the order
of poles and the global procedures coordinating computations are exactly
the same as for type $B_n$.

The only serious adaptations to the case of $D_n$ appears in the
procedure ${\tt CheckDvector(n,v)}$.
In fact now we check that

\begin{eqnarray*}
v_1+\cdots+v_i        &\geq&0\quad\mbox{for }1\leq i\leq n-1,\\
v_1+\cdots+v_{n-1}+v_n&\geq&0\quad\mbox{and is even},\\
v_1+\cdots+v_{n-1}-v_n&\geq&0\quad\mbox{and is even}.
\end{eqnarray*}

Other modifications are in procedures that are parent of
${\tt CheckDvector}$.
For example the procedure ${\tt CheckDwall}$ works exactly as
${\tt CheckBwall}$, but now calls ${\tt CheckDvector}$
instead of ${\tt CheckBvector}$.
See Section~\ref{sect.B.mpns}.


\section{Performance of the programs}
\label{sect.perf}

In this Section, we describe several tests of our programs
implementing the above MNS algorithms for types $A_n$, $B_n$,
$C_n$, $D_n$. The algorithm implementation is made with Maple.
Our programs are freely available at {\tt www.math.polytechnique.fr/$\sim$vergne/work/IntegralPoints.html}.
We  compare our results with 
the ones obtained by two previous algorithms:

\begin{itemize}
\item The Sp (for \emph{special permutations}) algorithm by
  Baldoni-DeLoera-Vergne~\cite{BalDeLoeVer03}, only for $A_n$.
\item The implementation ${\tt LattE}$ of Barvinok's
  algorithm~\cite{latte}, for every classical algebra;
\end{itemize}

These 
two methods also helped us to test our algorithms on various
examples.

\begin{figure}[ht]
\begin{center}
\includegraphics[height=6cm]{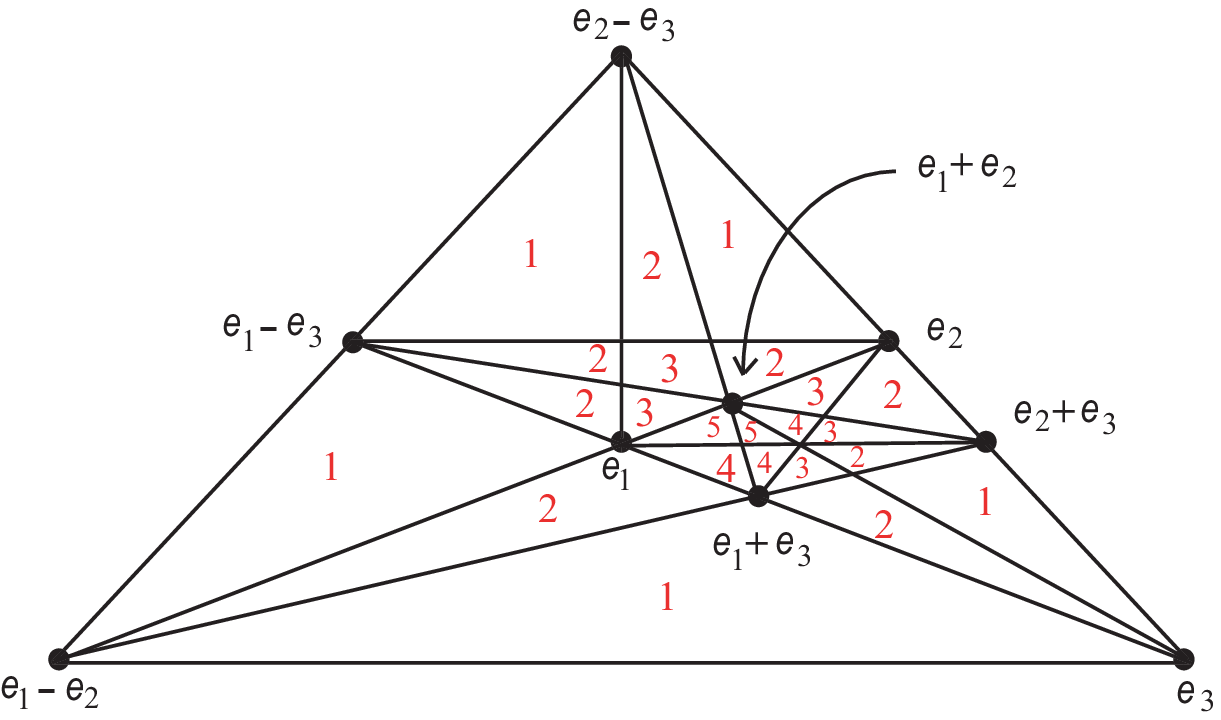}
\caption{Number of MNS containing any vector in a given chamber for
  $B_3$}
\label{figu.B3.MNS}
\end{center}
\end{figure}

Note that for our programs most of computation time is spent while
computing iterated residues. Indeed MNS computation is fast and
efficient. Note also that most of memory used by our programs
serves to store all fractions that occur in the iterated residue
process. The number of MNSs has a great influence on computation
time, since we sum over all MNSs. In any case it seems that the
deeper a vector is in the cone generated by positive roots, the
higher the number of MNSs is. This is morally bound to the fact
that there are more simplicial cones that might contain the
vector. In Figure~\ref{figu.B3.MNS}, we attach to every chamber
$\gc$ for $B_3$ the number of MNSs associated to any vector
$v\in\gc$.

Recall that the Sp algorithm relies on sums over a set
$\spa$ of \emph{special permutations} for a vector $a$.
The main advantage of our algorithms is that we compute fewer
iterated residues.
In fact the number of MNS seems to be smaller than the number of
special permutations that occur, for a given generic example.
But, examples at the end of Table~\ref{figu.A} show that a number of
MNSs considerably smaller than those of Sp's doesn't lead to a better
performance in time computation, even in the extreme case of just one
MNS.
Indeed this one residue computation can be  very time consuming due to
the substitutions $z_i=z_j$, which takes more time that the
substitutions $z_i=0$ used in the $\spa$ algorithm.
In the near future, we will improve this minor point. The MNSs method
should be better and is better in general.

During comparative tests, we figured out that one example
in~\cite{BalDeLoeVer03} has not been correctly copied from
draft.
More precisely in their Table~2 for complete graph $K_n$,
for the vector

$$a=(82275,33212,91868,-57457,47254,-64616,94854,-227390)$$
in the root lattice for $A_7$, the correct Kostant partition
number is the 103-digits integer

\begin{eqnarray*}
&&226040494681135377722281761934040091356424181\\
&&\quad\quad
  \quad\quad 242669497614801846058092972975120580334961426497
\end{eqnarray*}
and not only the first line of 45 digits.
The Kostant number and Ehrhart polynomials for this $a$ were
computed on a $1$GHz computer in $2,14\sec$ and $18,54\sec$
respectively, using $26$ special permutations.
Now with our programs running on a $1,13$GHz computer these times
drop to $1,38\sec$ and $2,50\sec$ respectively, using $14$ maximal
nested sets.
Similarly for the biggest example examined in~\cite{BalDeLoeVer03},
that is for the vector

\begin{eqnarray*}
&&a=(46398,36794,92409,-16156,29524,-68385,\\
&&\quad\quad\quad\quad 93335,50738,75167,-54015,-285809)
\end{eqnarray*}
in the root lattice for $A_{10}$, the 189-digits answer was
obtained in $2193\sec$ using $322$ special permutations,
whereas now we get the same result in $308\sec$ using $109$
maximal nested sets.

Table~\ref{figu.A} contains respective performances for $A_n$ of
${\tt LattE}$, Sp algorithms and our programs, a part the last four
examples that compare only the last program with ours. 
Tables~\ref{figu.B}--\ref{figu.D} contain respective performances
for $B_n$, $C_n$ and $D_n$ of  ${\tt LattE}$
and our programs.
We also indicated the number of special permutations (Sp) and
maximal nested sets (MNS).

Tests were performed on Pentium IV 1,13GHz computers with 1500 or
2000 mega-octets ($\mo$) of RAM memory.
We stopped several computations with {\tt LattE} when we figured out
that they would overcome computers' memory or take too much time with
respects to the other algorithms; in this case we indicate the time
spent and the number of mega-octets used by the computer.

\begin{figure}[ht]
\begin{center}
\rotatebox[]{90}{
$\begin{array}{||c||c|c|c||} \hline \hline
\mbox{Root lattice element}      &\mbox{{\tt LattE}}&\mbox{Sp}&\mbox{MNS}\\ \hline \hline
(2215,571,4553,-600,-6739)                           &    1,6\sec&    <0,1\sec,  \mbox{4 Sp}&      <0,1\sec, \mbox{3 MNS}\\ \hline
(6440,-4866,6174,-5683,7112,-9177)                   &    2,0\sec&    <0,1\sec,  \mbox{4 Sp}&       0,1\sec, \mbox{1 MNS}\\ \hline
(5067,3639,-3103,435,-729,2267,-7576)                &   61,5\sec&     0,1\sec, \mbox{12 Sp}&       0,3\sec, \mbox{8 MNS}\\ \hline
(2232,-1656,7452,99,601,-2870,-2908,-2950)           &  808,8\sec&     1,6\sec, \mbox{56 Sp}&       1,2\sec, \mbox{9 MNS}\\ \hline
(4060,183,-4211,5914,2790,-5360,-1730,3916,-5562)    & 1646,2\sec&     4,8\sec, \mbox{40 Sp}&       0,3\sec, \mbox{2 MNS}\\ \hline
(4058,-1343,-2236,7114,1909,                         & \mbox{--} &    42,5\sec &       2,3\sec\\
\quad\quad\quad\quad\quad -5696,193,5298,-689,-8608) &           & \mbox{64 Sp}&  \mbox{8 MNS}\\ \hline
(1388,4024,-1586,-1135,5998,-6067,                   & \mbox{--} &  1162,9\sec &      12,1\sec\\
\quad\quad\quad 3562,-4599,7818,-2542,-6861)         &           &\mbox{256 Sp}&  \mbox{6 MNS}\\ \hline
(1094,-11,-75,1,-1009)                               &    0,6\sec&    <0,1\sec,  \mbox{4 Sp}&      <0,1\sec, \mbox{1 MNS}\\ \hline
(1034,49,-75,25,-33,-1000)                           &    7,1\sec&    <0,1\sec, \mbox{16 Sp}&       0,3\sec, \mbox{6 MNS}\\ \hline
(1022,36,33,-53,-21,-1,-1016)                        &  182,1\sec&     0,3\sec, \mbox{40 Sp}&       1,4\sec, \mbox{20 MNS}\\ \hline
(1099,-99,77,-15,-29,24,36,-1093)                    &  337,0\sec&     0,3\sec,  \mbox{8 Sp}&       0,3\sec, \mbox{4 MNS}\\ \hline
(1050,-36,5,-130,-16,43,20,91,-1027)                 & 3764,1\sec&     1,6\sec, \mbox{20 Sp}&       0,7\sec, \mbox{3 MNS}\\ \hline
(1079,-64,28,11,-48,5,-4,25,20,-1052)                & \mbox{--} &    23,8\sec, \mbox{40 Sp}&       5,0\sec, \mbox{12 MNS}\\ \hline
(1052,-46,-52,25,-21,69,-26,25,-43,24,-1007)         & \mbox{--} &   896,4\sec, \mbox{216 Sp}&      41,6\sec,\mbox{32 MNS}\\ \hline
(31011,1000,600,500,-500,-600,-1000,-31011)          &12832,8\sec&     3,1\sec &      18,0\sec\\
                                                     &   1500\mo &\mbox{206 Sp}&\mbox{137 MNS}\\ \hline
(31011,10000,6000,5000,0,-5000,-6000,-10000,-31011)  & >23000\sec&    60,4\sec &    1865,8\sec\\
                                                     &           &\mbox{898 Sp}&\mbox{548 MNS}\\ \hline
(46398,36794,92409,-16156,29524,-68385,              & \mbox{--} &  2193,2\sec &     308,5\sec\\
\quad\quad\quad 93335,50738,75167,-54015,-285809)    &          &\mbox{322 Sp}&\mbox{109 MNS}\\ \hline
(37,-9,-7,-6,-5,-4,-3,-2,-1)                         & >12000\sec&     7,0\sec &     213,6\sec\\
                                                     &  >2400\mo &\mbox{128 Sp}&  \mbox{1 MNS}\\ \hline
\end{array}$
}
\caption{Computation time for 
${\tt LattE}$,
  B-DL-V and our programs, for $A_n$}
\label{figu.A}
\end{center}
\end{figure}

\begin{figure}[ht]
\begin{center}
$\begin{array}{||c||c|c||} \hline \hline
\mbox{Root lattice element}      &\mbox{{\tt LattE}}&       \mbox{MNS} \\ \hline \hline
(1388,4024,3826)              &    0,8\sec&      <0,1\sec\\
                              &           &  \mbox{3 MNS}\\ \hline
(2691,5998,-6067,6184)        &    2,6\sec&       0,1\sec\\
                              &           &  \mbox{1 MNS}\\ \hline
(1585,7818,-2542,-2803,2715)  &  214,9\sec&       3,0\sec\\
                              &           &  \mbox{2 MNS}\\ \hline
(479,7114,1909,-5696,193,9297)&16369,6\sec&      27,5\sec\\
                              &           &  \mbox{8 MNS}\\ \hline
(1070,1006,-37)               &    0,9\sec&       0,1\sec\\
                              &           &  \mbox{3 MNS}\\ \hline
(1082,947,27,42)              &   22,9\sec&       1,2\sec\\
                              &           & \mbox{15 MNS}\\ \hline
(1047,974,20,44,-35)          & 1939,9\sec&      21,7\sec\\
                              &           & \mbox{51 MNS}\\ \hline
(1015,1082,-37,-21,-28,14)    &  >7000\sec&     378,0\sec\\
                              &   >1500\mo& \mbox{26 MNS}\\ \hline
\end{array}$
\caption{Computation time for 
${\tt LattE}$
  and our programs, for $B_n$}
\label{figu.B}
\end{center}
\end{figure}

\begin{figure}[ht]
\begin{center}
$\begin{array}{||c||c|c||} \hline \hline
\mbox{Root lattice element}      &\mbox{{\tt LattE}}&        \mbox{MNS}\\ \hline \hline
(1388,4024,7652)               &   0,8\sec&      <0,1\sec\\
                               &          &  \mbox{1 MNS}\\ \hline
(2691,5998,-6067,12368)        &   2,8\sec&       0,1\sec\\
                               &          &  \mbox{1 MNS}\\ \hline
(1585,7818,-2542,-2803,5430)   & 163,0\sec&       1,4\sec\\
                               &          &  \mbox{1 MNS}\\ \hline
(479,7114,1909,-5696,192,18594)& >5400\sec&      65,3\sec\\
                               &   >900\mo&  \mbox{8 MNS}\\ \hline
(1038,22,-2)                   &   0,8\sec&      0,1\sec\\
                               &          & \mbox{3 MNS}\\ \hline
(1021,37,-40,178)              &  12,2\sec&      0,5\sec\\
                               &          & \mbox{4 MNS}\\ \hline
(1051,-45,26,-5,-131)          & 195,4\sec&      2,8\sec\\
                               &          & \mbox{6 MNS}\\ \hline
(1024,6,60,-6,-42,52)          &>10800\sec&    1292,4\sec\\
                               &  >2000\mo&\mbox{42 MNS}\\ \hline
\end{array}$
\caption{Computation time for 
${\tt LattE}$
  and our programs, for $C_n$}
\label{figu.C}
\end{center}
\end{figure}

\begin{figure}[ht]
\begin{center}
$\begin{array}{||c||c|c||} \hline \hline
\mbox{Root lattice element}      &\mbox{{\tt LattE}}&      \mbox{MNS}\\ \hline \hline
(8608,-305,183)                     &   0,3\sec &      <0,1\sec\\
                                    &           &  \mbox{1 MNS}\\ \hline
(32,5914,6166,-5360)                &   1,5\sec &      <0,1\sec\\
                                    &           &  \mbox{1 MNS}\\ \hline
(1646,3916,-3330,6372,7452)         &  18,0\sec &       0,5\sec\\
                                    &           &  \mbox{2 MNS}\\ \hline
(8127,601,-2870,-2908,10823,3639)   & 313,5\sec &       3,1\sec\\
                                    &           &  \mbox{2 MNS}\\ \hline
(1009,1106,-9)                      &   0,2\sec &      <0,1\sec\\
                                    &           &  \mbox{1 MNS}\\ \hline
(1074,959,64,77)                    &   3,0\sec &       0,3\sec\\
                                    &           &  \mbox{6 MNS}\\ \hline
(1100,973,2,-1,-60)                 & 100,2\sec &       3,1\sec\\
                                    &           & \mbox{18 MNS}\\ \hline
(1096,965,-54,68,-34,-1)            &7076,7\sec &     763,3\sec\\
                                    &           & \mbox{47 MNS}\\ \hline
\end{array}$
\caption{Computation time for  ${\tt LattE}$
  and our programs, for $D_n$}
\label{figu.D}
\end{center}
\end{figure}


\end{document}